\documentclass[letter,twocolumn]{autart}    

\newcommand{\new}[1]{{\color{blue} #1}}
\renewcommand{\new}[1]{#1}

\usepackage[normalem]{ulem}
\usepackage{graphics} 
\usepackage{epsfig} 
\usepackage{amsmath} 
\usepackage{amssymb}  
\usepackage{cite}
\usepackage{color}
\usepackage{graphicx}
\usepackage{algpseudocode}
\usepackage{algorithm}
\usepackage[numbers]{natbib}
\usepackage{wrapfig} 

\newtheorem{theorem}{Theorem}[section]
\newtheorem{definition}[theorem]{Definition}
\newtheorem{assumption}[theorem]{Assumption}
\newtheorem{lemma}[theorem]{Lemma}

\newtheorem{remark}[theorem]{Remark}
\newtheorem{example}[theorem]{Example}

\newtheorem{proposition}[theorem]{Proposition}  
\newtheorem{problem}[theorem]{Problem}

\usepackage{enumitem}
\renewcommand{\theenumi}{\alph{enumi}}

\newcommand{\naturals}{\mathbb{N}}
\newcommand{\real}{\mathbb{R}}
\newcommand{\cplx}{\mathbb{C}}

\newcommand{\range}{\mathcal{R}}

\newcommand{\Ac}{\mathcal{A}}

\newcommand{\Fc}{\mathcal{F}}
\newcommand{\Gc}{\mathcal{G}}
\newcommand{\Hc}{\mathcal{H}}

\newcommand{\Kc}{\mathcal{K}}
\newcommand{\Lc}{\mathcal{L}}
\newcommand{\Mc}{\mathcal{M}}

\newcommand{\Pc}{\mathcal{P}}

\newcommand{\Sc}{\mathcal{S}}
\newcommand{\Tc}{\mathcal{T}}
\newcommand{\Uc}{\mathcal{U}}
\newcommand{\Vc}{\mathcal{V}}

\newcommand{\Xc}{\mathcal{X}}

\newcommand{\Pf}{\mathfrak{P}}

\newcommand{\Kapprox}{K_{\operatorname{approx}}}
\newcommand{\Kedmd}{K_{\operatorname{EDMD}}}

\newcommand{\until}[1]{\{1,\dots,#1\}}

\newcommand{\psix}{\Psi(X)}
\newcommand{\psixp}{\Psi(X^+)}

\newcommand{\Span}{\operatorname{span}}

\newcommand{\Tr}{\operatorname{Tr}}

\newcommand{\dist}{{\operatorname{dist}}}

\newcommand{\edmd}{\operatorname{EDMD}}

\newcommand{\ic}{\mathcal{I}_C}

\newcommand{\restr}[2]{#1 \!\! \restriction_{#2}}

\newcommand{\aug}{\operatorname{aug}}
\newcommand{\Faug}{\Fc^{\aug}}
\newcommand{\Taug}{\Tc^{\aug}}
\newcommand{\Kaug}{\Kc^{\aug}}
\newcommand{\innerprod}[2]{\langle #1, #2 \rangle}

\newcommand{\longthmtitle}[1]{\mbox{}{\textit{(#1):}}}
\newcommand{\setdef}[2]{\{#1 \; | \; #2\}}

\newcommand{\oprocendsymbol}{\hbox{$\square$}}
\newcommand{\oprocend}{\relax\ifmmode\else\unskip\hfill\fi\oprocendsymbol}
\def\eqoprocend{\tag*{$\square$}}

\allowdisplaybreaks

\graphicspath{{epsfiles/}}

\begin{document}
\begin{frontmatter}

\title{Modeling Nonlinear Control Systems via Koopman Control Family: 
Universal Forms and Subspace Invariance Proximity}

\thanks{This work was supported by ONR Award N00014-23-1-2353 and
NSF Award IIS-2007141.}

\vspace{-10pt}

\author[First]{Masih Haseli}%
\author[First]{\quad Jorge Cort\'es}

\address[First]{Department of Mechanical and Aerospace
Engineering, University of California, San Diego,
\{mhaseli,cortes\}@ucsd.edu}

\begin{abstract}
This paper introduces the Koopman Control Family (KCF), a
mathematical framework for modeling general (not necessarily control-affine) discrete-time nonlinear
control systems with the aim of providing a solid theoretical
foundation for the use of Koopman-based methods in systems with
inputs.  We demonstrate that the concept of KCF captures the
behavior of nonlinear control systems on a (potentially
infinite-dimensional) function space. By employing a generalized
notion of subspace invariance under the KCF, we establish a
universal form for finite-dimensional models, which encompasses the
commonly used linear, bilinear, and linear switched models as
specific instances. In cases where the subspace is not invariant
under the KCF, we propose a method for approximating models in
general form and characterize the model's accuracy using the concept
of invariance proximity.  We end by discussing how
the proposed framework naturally lends itself to data-driven
modeling of control systems.
\end{abstract}


\end{frontmatter}

\section{Introduction}
The Koopman operator approach to dynamical systems has gained
widespread attention in recent years. While traditional state-space
methods for nonlinear systems rely on the description of system
trajectories, the Koopman viewpoint provides an equivalent formulation
of the system behavior using a linear operator acting on a vector
space of functions.
The Koopman operator framework yields beneficial algebraic
constructions that can be leveraged for efficient computational
learning and prediction.  These benefits have motivated researchers to
consider extending the framework to control systems.  However, unlike
the case of \new{systems without input,}
this has proven difficult due to the fact that the input's role is
fundamentally different from the state's role: without a priori
knowledge of the input signal, there is not enough information to
predict the system trajectories since the choice of input can
drastically alter system behavior.
Our aim is to provide a comprehensive mathematical framework for
Koopman operator-based modeling of control systems.

\textit{Literature Review:} The Koopman operator~\citep{BOK:31}
provides an alternative description of nonlinear 
systems \new{(without input)} that encodes the system behavior through the evolution of
functions (a.k.a., observables) belonging to a vector space. Even
though the system might be nonlinear, the Koopman operator is always
linear. This inherent linearity gives rise to favorable algebraic
properties, leading to powerful tools to analyze complex dynamical
systems~\citep{IM:05,CWR-IM-SB-PS-DSH:09,MB-RM-IM:12} for which typical
state-space and geometric methods are cumbersome. However, linearity
comes at the expense of the infinite--dimensional nature of the
operator. To make possible its direct and efficient implementation on
digital computers, one needs to develop finite-dimensional
descriptions for it. This generally relies on the concept of subspace
invariance~\citep{SLB-BWB-JLP-JNK:16}. If a finite-dimensional subspace
is invariant under the operator, then one can restrict the operator to
the subspace and represent its action with a matrix given a chosen
basis. This has led to a search for invariant subspaces through a
variety of techniques, including the direct identification of
eigenfunctions (which span invariant
subspaces)~\citep{MK-IM:20,EK-JNK-SLB:21}, optimization and neural
network-based
methods~\citep{NT-YK-TY:17,BL-JNK-SLB:18,EY-SK-NH:19,SEO-CWR:19,SP-NAM-KD:21,MS:21-L4DC,MK:23-tac},
and efficient algebraic searches~\citep{MH-JC:22-tac,MH-JC:21-tcns}.

Even without a finite-dimensional invariant subspace available, one
can still approximate the action of the Koopman operator on any
finite-dimensional subspace via an orthogonal projection.  Prominent
data-driven methods in this category are Dynamic Mode Decomposition
(DMD)~\citep{PJS:10,JHT-CWR-DML-SLB-JNK:13} and its variant, Extended
Dynamic Mode Decomposition
(EDMD)~\citep{MOW-IGK-CWR:15,MK-IM:18}. In
addition,~\citep{SM-NC-IM:18} provides methods to handle
time-varying systems.
For such methods, criteria that must be balanced to choose the
finite-dimensional space are the relevance of the dynamical
information captured by the subspace and the accuracy of the resulting
approximation.  The work~\citep{MH-JC:23-csl} provides a tight upper
bound on the worst-case prediction error on a subspace, providing a
measure of the quality of the subspace independently of the chosen
basis.  The works~\citep{HL-DMT:20,FN-SP-FP-MS-KW:23} provide several
error bounds for accuracy of DMD, EDMD, and extensions to
Koopman-based control models.
The work~\citep{MH-JC:23-auto} provides an algebraic algorithm to
approximate Koopman-invariant subspaces of an arbitrary
finite-dimensional space with tunable accuracy.

The algebraic properties of the Koopman operator have been used in a
myriad applications, including fluid
dynamics~\citep{CWR-IM-SB-PS-DSH:09}, stability
analysis~\citep{AM-IM:16,BY-IRM:21,SAD-AMV-CJT:22,CMZ-AM:21}, reachability
analysis~\citep{SB-SB-PSD-ARG-KP:21,BU-DT-UV:22,HB-AA-ZL-NO:23,WS-NS-MK-YTC-SH:23},
safety-critical control~\citep{CF-YC-ADM-JWB:20,VZ-EB:23,MB-DP:22},
and robotics~\citep{DB-XF-RBG-CDR-RC:20,LS-KK:21,GM-IA-TDM:23}.
The fact that the Koopman operator is only formally defined for
 systems \new{without input} has not been an obstacle for the development of
many data-based methods inspired by it to construct low-complexity
representations of control systems.  Many advances do not directly
require an operator-theoretic approach for the open control system,
but instead rely on the Koopman operator of the unforced system and
address control design based on control Lyapunov
functions~\citep{AN-SHS-JSK:22,VZ-EB:23-neurocomputing} or feedback
linearization~\citep{DG-VK-FP:24}.  A significant amount of attention
has been devoted to deriving finite-dimensional forms by lifting to
higher dimensions.  Due to their simplicity, lifted linear models are
the most popular in the
literature~\citep{SLB-BWB-JLP-JNK:16,MK-IM-automatica:18}
and lead to highly efficient computational algorithms. Such models,
however, are not capable of capturing certain structures, such as
terms containing the products of inputs and states, which are
prevalent in control-affine nonlinear systems. For these, the
works~\citep{BH-XM-UV:18,SP-SEO-CWR:20,DG-DAP:21,RS-JB-FA:23} propose
the use of bilinear lifted models based on geometric arguments relying
on the control-affine structure. The work~\citep{LCI-RT-MS:22} proposes
a different lifted form based on invariant subspaces for the Koopman
operator associated with the unforced system. An interesting
alternative is to model the system by switching between several
Koopman operators, each associated with the system under a different
constant~\citep{SP-SK:19,MB-JPH:23} or piecewise
constant~\citep{AS-DE:17,AS-AM-DE:18} input signal.  The
work~\citep{MK-IM-automatica:18} takes a different approach and
considers the system behavior under all possible infinite input
sequences, defining a Koopman operator for the control system on a
function space whose members' domain is the Cartesian product of the
state space and all possible input sequences. This is perhaps the most
direct approach in terms of an operator-theoretic viewpoint for
controls systems. However, given the reliance on infinite input
sequences, working with finite-time restrictions should be done with
care. Here, we take a different operator-theoretic approach to capture
the behavior of control systems that we find easier to work with on
finite-dimensional subspaces with only finitely many input steps
available.

\textit{Statement of Contributions:} Our goal here is to provide a
solid theoretical framework to model general discrete-time nonlinear
\emph{control} systems based on Koopman operator theory.  The starting
point of our approach is the observation made in the literature that
if the input is a constant signal, then the control system becomes 
\new{a system without input.}
Therefore, any Koopman-theoretic model for the control system must
reduce to the conventional Koopman operator associated with the
corresponding 
\new{system without input.}
 Motivated by this, we define the
concept of \emph{Koopman control family (KFC)} as the collection of
all Koopman operators associated with constant-input 
dynamics derived from the control system. We show that the KCF can
completely capture the control system behavior on a potentially
infinite-dimensional function space.  Since dealing with
infinite-dimensional operators is computationally intractable, we
provide a theoretical structure for finite-dimensional models whose
construction is based on projection operators, analogous to the case
of \new{systems without input.}
To find a general finite-dimensional form for
Koopman-based models for the control system, we rely on a generalized
notion of subspace invariance. Specifically, we show that on a
common-invariant subspace for the KCF, the finite-dimensional model
always has a specific \emph{``input-state separable''}
form. Remarkably, the linear, bilinear, and switched linear models
commonly used in the literature are all special cases of the
input-state separable form.  Since KCF contains uncountably many
operators (given the infinite choices for the constant input signal),
directly finding a common invariant subspace is challenging.  To
tackle this, we parametrize the KCF with one operator, termed
augmented Koopman operator, and show that invariant subspaces under
this augmented operator lead to common invariant subspaces for the
KCF. As a result, the problem of working with uncountably many
operators simplifies to working with a single linear operator on an
extended function space.  Similarly to the case of \new{systems without input,}
finding an exact and informative finite-dimensional common invariant
subspace for the KCF is generally challenging and in some cases might
not even exist. To address this, we define the concept of invariance
proximity under an operator, which enables us to extended our results
to approximations on non-invariant subspaces and provide a bound on
the accuracy of the resulting approximated models for the control
system.  Our final contribution shows how the results can directly be
used in data-driven modeling applications.

\section{Preliminaries}\label{sec:preliminaries}
In this section\footnote{The symbols $\naturals$, $\real$, and
  $\cplx$, represent the sets of natural, real, and complex numbers,
  resp.  Given $A \in \cplx^{m \times n}$, we denote its transpose,
  pseudo-inverse, conjugate transpose, Frobenius norm and range space
  by $A^T$, $A^\dagger$, $A^H$, $\|A\|_F$, and $\range(A)$, resp. If
  $A$ is square, $A^{-1}$ and $\Tr(A)$ denote its inverse and trace,
  resp. When all eigenvalues of $A$ are real, $\lambda_{\min}(A)$ and
  $\lambda_{\max}(A)$ denote the smallest and largest eigenvalues
  of~$A$.  We use $I_m$ and $\mathbf{0}_{m\times n}$ to denote the
  $m \times m$ identity matrix and $m\times n$ zero matrix (we drop
  the indices when appropriate). We denote the $2$-norm of the vector
  $v \in \cplx^n$ by $\|v\|_2$. Given sets $S_1$ and $S_2$, their
  union and intersection are represented by $S_1 \cup S_2$ and
  $S_1 \cap S_2$, resp. Also, $S_1 \subseteq S_2$ and
  $S_1 \subsetneq S_2$ mean that $S_1$ is a subset and proper subset
  of $S_2$, resp. Given the vector space $\Vc$ defined on the field
  $\cplx$, $\dim \Vc$ denotes its dimension. Given a set
  $\Sc \subseteq \Vc$, $\Span(\Sc)$ is a vector space comprised of all
  linear combinations of elements in $\Sc$. Given vector spaces
  $\Vc_1$ and $\Vc_2$, we define
  $\Vc_1+\Vc_2 := \setdef{v_1 + v_2}{v_1 \in \Vc_1, \, v_2 \in
    \Vc_2}$.  Given functions $f$ and $g$ with appropriate domains and
  co-domains, $f \circ g$ denotes their composition.  Let
  $f: A \times B \to C$ be a multivariable function yielding $f(a,b)$
  for $(a,b) \in A\times B$. Then, $\restr{f}{b=b^*}:A \to C$ is
  defined as $\restr{f}{b=b^*}(a) := f(a,b^*)$ for $a \in A$. If $F$
  consists of multivariable functions of the form
  $f: A \times B \to C$, then
  $\restr{F}{b = b^*}:=\setdef{\restr{f}{b=b^*}}{f \in F}$. Given a
  positive measure $\mu$ on a set $A$ and functions
  $f,g: A \to \cplx$, we define their $L_2$ inner product as
  $\innerprod{f}{g}_{L_2(\mu)}:= \int_A f(x)\bar{g}(x) d\mu(x)$, where
  $\bar{g}$ is the complex conjugate of $g$. The $L_2$ norm of $f$ is
  defined as $\|f\|_{L_2(\mu)} = \sqrt{\innerprod{f}{f}_{L_2(\mu)}}$.
  We drop the dependency on the measure $\mu$ when the context is
  clear.}, we review notions and results regarding the Koopman
operator, Extended Dynamic Mode Decomposition, and the concept of
consistency index.

\subsection{Koopman Operator}\label{sec:prelim-Koopman}
Here, we briefly explain the Koopman operator associated with a
dynamical system and its properties following the terminology
in~\cite{MB-RM-IM:12}. Consider a discrete-time system
\begin{align}\label{eq:dynamical-sys}
x^+ = T(x), 
\end{align}
with state space $\Xc \subseteq \real^n$, where $T: \Xc \to \Xc$
is the function describing the system behavior.  Consider a
linear function space $\Fc$ (defined on the field $\cplx$) comprised
of functions form $\Xc$ to $\cplx$ and assume it is closed under
composition with $T$, i.e., $f \circ T \in \Fc$ for all $f \in
\Fc$. We define the Koopman operator $\Kc: \Fc \to \Fc$ as
\begin{align}\label{eq:Koopman-def}
\Kc f = f \circ T.
\end{align}
It is easy to verify that~\eqref{eq:Koopman-def} is a linear operator,
i.e.,
\begin{align}\label{eq:Koopman-spatial-linear}
\Kc (\alpha f + \beta g) = \alpha \Kc f + \beta \Kc g, \; \forall
f,g \in \Fc, \; \forall \alpha, \beta \in \cplx. 
\end{align}
The Koopman operator's action on a given function can be viewed as
pushing forward the function values across all system trajectories
by one time step. We can repeatedly apply the Koopman operator on a
function $f \in \Fc$ to predict its evolution on any system trajectory
$\{x(i)\}_{i=0}^\infty$ as
\begin{align}\label{eq:multistep-Koopman-nocontrol}
f(x(i)) = \Kc^i f(x(0)), \quad \forall i \in \naturals_0.
\end{align}
Since $\Kc$ is a linear operator, we can define its
eigendecomposition. We say the function $\phi \in \Fc$ is an
\emph{eigenfunction} of $\Kc$ with \emph{eigenvalue} $\lambda$ if
\begin{align}\label{eq:Koopman-eig-def}
\Kc \phi = \lambda \phi.
\end{align}
By comparing~\eqref{eq:multistep-Koopman-nocontrol}
and~\eqref{eq:Koopman-eig-def}, one can see that Koopman
eigenfunctions evolve linearly on system trajectories,
\begin{align}\label{eq:eigs-temporal-linearity}
\phi(x(i)) = \lambda \phi(x(i-1)), \quad \forall i \in \naturals.
\end{align}
We refer to~\eqref{eq:eigs-temporal-linearity} as \emph{temporal
linear evolution} of eigenfunctions. This temporal linearity of
eigenfunctions combined with the
linearity~\eqref{eq:Koopman-spatial-linear} of the operator on the
space $\Fc$ enables us to linearly predict function values on system
trajectories. Specifically, given eigenfunctions
$\{\phi_k\}_{k=1}^{N_k}$ with corresponding eigenvalues
$\{\lambda_k\}_{k=1}^{N_k}$, one can write the evolution of the
function $f = \sum_{k=1}^{N_k} c_k \phi_k$ on the system
trajectories as
\begin{align*}
f(x(i)) = \sum_{k=1}^{N_k} c_k \lambda_k^i \phi_k (x(0)), \;
\forall i \in \naturals_0. 
\end{align*}
This equation is of utmost importance since it provides a linear
structure facilitating the prediction of nonlinear
systems~\cite{MK-IM-automatica:18,MK-IM:20} as well as learning the
system behavior from
data~\cite{CWR-IM-SB-PS-DSH:09,MH-JC:22-tac,MH-JC:23-auto}.
One should keep in mind that, in general, to capture the full state of
the system, one might need $\Fc$ to be infinite dimensional since it
must be closed under composition with~$T$.

Next, we define the concept of subspace invariance under the Koopman
operator. 
A subspace $\Gc \subseteq \Fc$ is \emph{Koopman invariant}
if $\Kc f \in \Gc$ for all $f \in \Gc$. Koopman eigenfunctions
trivially span invariant subspaces. 

\begin{remark}\longthmtitle{Simplifying Notation For Vector-Valued
Functions}\label{r:notation-overload}
{\rm 
For convenience, we introduce some notation simplifying the
operation of the Koopman operator on finite-dimensional spaces.  Let
$\Psi: \Xc \to \cplx^s$ be a vector-valued map represented as
$\Psi(\cdot) = [\psi_1(\cdot),\ldots, \psi_s(\cdot)]^T$, where
$\psi_i: \Xc \to \cplx$ for all $i \in \until{s}$. We define the
span of the elements of $\Psi$ and action of Koopman operator on the
elements of $\Psi$ as
\begin{align*}
\Span(\Psi) & := \Span(\{\psi_1,\ldots, \psi_s\}),
\\
\Kc \Psi & := [\Kc \psi_1, \ldots, \Kc \psi_s]^T = \Psi \circ T.
\end{align*}
Given a finite-dimensional subspace $\Hc \subset \Fc$, we
often describe a basis for it by a vector-valued map
$\Phi: \Xc \to \cplx^{\dim(\Hc)}$ satisfying
$\Hc = \Span(\Phi)$.  \oprocend }
\end{remark}

An important property of finite-dimensional Koopman-invariant
subspaces is that one can capture the action of the operator by a
matrix once a basis is chosen. Formally, given the invariant subspace
$\Sc \subseteq \Fc$ with basis $\Psi: \Xc \to \cplx^{\dim(\Sc)}$,
there exists a unique matrix
$K \in \cplx^{\dim(\Sc) \times \dim(\Sc)}$ such that
\begin{align}\label{eq:Koopman-matrix-invariant}
\Kc \Psi= \Psi \circ T= K \Psi.
\end{align}
This equation combined with the linearity of the operator allow us to
easily calculate the action of the operator on any function in
$\Sc$. Formally, given any function $f \in \Sc$ with description
$f (\cdot) = w^T \Psi(\cdot)$ where $w \in \cplx^{\dim(\Sc)}$, one has
\begin{align}\label{eq:Koopman-predictor-invariant}
\Kc f = w^T K \Psi.
\end{align}

The concept of subspace invariance is of utmost importance since it
allows us to operate on finite-dimensional subspaces and use numerical
matrix computation for prediction, as laid out in
equations~\eqref{eq:Koopman-matrix-invariant}-\eqref{eq:Koopman-predictor-invariant}.

Even if the subspace $\Sc \subseteq \Fc$ is \emph{not} invariant under
$\Kc$, it is still possible to use the notion of subspace
invariance to \emph{approximate} the action of $\Kc$ on $\Sc$. To
achieve this, one usually utilizes $\Pc_{\Sc}: \Fc \to \Fc$, the
orthogonal projection operator (given an inner product on $\Fc$) on
$\Sc$. Observe that the space $\Sc$ is invariant under the operator
$\Pc_{\Sc} \Kc: \Fc \to \Fc$; hence,
equations~\eqref{eq:Koopman-matrix-invariant}-\eqref{eq:Koopman-predictor-invariant}
are valid when we substitute in them the Koopman operator $\Kc$ by
$\Pc_{\Sc} \Kc$. Let $\Kapprox$ be the matrix calculated by applying
equation~\eqref{eq:Koopman-matrix-invariant} to the operator
$\Pc_{\Sc} \Kc$. Then, this matrix provides an approximation for the
action of $\Kc$ on $\Sc$ as follows
\begin{align}\label{eq:Koopman-matrix-noninvariant}
\Kc \Psi= \Psi \circ T \approx \Pc_{\Sc} \Kc \Psi= \Kapprox \Psi.
\end{align}
Moreover, the analogous but approximated version of function
prediction in~\eqref{eq:Koopman-predictor-invariant} is given by
\begin{align}\label{eq:Koopman-predictor-noninvariant}
\Kc f \approx \Pc_{\Sc} \Kc f = w^T \Kapprox \Psi.
\end{align}

\begin{remark}\longthmtitle{General Linear Form and Subspace
Invariance}\label{r:invariance-general-form}
{\rm When dealing with the action of the Koopman operator on
  finite-dimensional spaces, we use linear models that are either
  exact,
  cf.~\eqref{eq:Koopman-matrix-invariant}-\eqref{eq:Koopman-predictor-invariant},
  or approximated,
  cf.~\eqref{eq:Koopman-matrix-noninvariant}-\eqref{eq:Koopman-predictor-noninvariant}.
  In either case the model has the same form. It is in this sense we
  say the general finite-dimensional form of Koopman-based models is
  linear. Note that this general form is a consequence of the notion
  of subspace invariance (whether the subspace is actually
  Koopman-invariant or not). \oprocend }
\end{remark}

\subsection{Extended Dynamic Mode Decomposition}\label{subsec:EDMD}
In many engineering applications, the system dynamics is unknown and
we only have access to data from the system trajectories.  The
Extended Dynamic Mode Decomposition (EDMD)
method~\cite{MOW-IGK-CWR:15} reviewed here
uses data to approximate the action of the Koopman operator on a given
\emph{finite-dimensional} space of functions.

\begin{remark}\longthmtitle{Use of Real-valued Basis Functions in
Data-Driven Applications}\label{r:real-dictionary}
{\rm All the
systems in this paper are defined on state and input spaces with
real-valued elements. Consequently, the Koopman operator's action on
any pair of complex-conjugate functions leads to another
complex-conjugate pair, which can be captured by a pair of
real-valued functions. Hence, even though we develop our theory
based on complex functions, in data-driven applications, we work
with bases with real-valued elements to simplify the numerical
operations, without loss of generality\footnote{Given a
vector-valued function $\Psi$ with real-valued elements,
$\Span(\Psi)$ contains complex-valued functions since we employ
$\cplx$ as the underlying field.}.  \oprocend
}
\end{remark}

To specify the function space, EDMD uses a dictionary comprised of
$N_{\Psi}$ functions form $\Xc$ to $\real$. Formally, we define our
dictionary as a vector-valued function
\begin{align*}
\Psi(\cdot) = [\psi_1(\cdot),\ldots, \psi_{N_{\Psi}}(\cdot)]^T,
\end{align*}
where $\psi_1, \ldots, \psi_{N_{\Psi}} \in \Fc$ are the dictionary
elements. To approximate the behavior of the Koopman operator on
$\Span(\Psi)$, EDMD uses data snapshots from the trajectories in two
matrices $X,X^+ \in \real^{n \times N}$ such that
\begin{align}\label{eq:data-snapshots}
x_i^+ = T(x_i), \; \forall i \in \until{N},
\end{align}
where $x_i$ and $x_i^+$ are the $i$th columns of matrices $X$ and
$X^+$, resp. For convenience, we define the action of the
dictionary on data matrix $X$ (similarly for any data matrix)~as
\begin{align*}
  \psix = [\Psi(x_1), \Psi(x_2), \ldots, \Psi(x_n)] \in \real^{N_{\psi} \times N}.
\end{align*}
Note that based on~\eqref{eq:Koopman-def}
and~\eqref{eq:data-snapshots}, one can see
$ \psixp = \Psi \circ T(X) = \Kc \Psi(X)$.  Hence, the dictionary
matrices $\psix$ and $\psixp$ capture the behavior of the Koopman
operator on $\Span(\Psi)$. EDMD approximates the action of the
operator by solving a least-squares problem
\begin{align}\label{eq:EDMD-optimization}
\underset{K}{\text{minimize}} \| \psixp -  K \psix \|_F,
\end{align} 
with the following closed-form solution
\begin{align}\label{eq:EDMD-closed-form}
\Kedmd = \psixp  \psix^\dagger.
\end{align}
Throughout this paper, we make the following assumption.

\begin{assumption}\longthmtitle{Full Rank Dictionary
Matrices}\label{a:full-rank}
$\psix$ and $\psixp$ have full row rank.  \oprocend
\end{assumption}

Note that Assumption~\ref{a:full-rank} implies that the element of
$\Psi$ are linearly independent. Also, it implies that data in $X$ and
$X^+$ are diverse enough to distinguish between functions in
$\Span(\Psi)$. Moreover, if Assumption~\ref{a:full-rank} holds,
$\Kedmd$ is the unique solution of~\eqref{eq:EDMD-optimization}.

The matrix $\Kedmd$ captures relevant information about the system
behavior and can be used to approximate the action of the Koopman
operator on $\Span(\Psi)$. Formally, we define the EDMD predictor of
$\Kc \Psi$ by EDMD as
\begin{align}\label{eq:EDMD-dictionary-predictor}
\Pf_{\Kc \Psi}^{\edmd} := \Kedmd \Psi.
\end{align}
Similarly, for an arbitrary function $f \in \Span(\Psi)$ with
description $f(\cdot) = w^T \Psi(\cdot)$ for $w \in \cplx^{N_{\Psi}}$,
one can define the EDMD predictor of $\Kc f$ as
\begin{align}\label{eq:EDMD-function-predictor}
\Pf_{\Kc f}^{\edmd} := w^T \Kedmd \Psi.
\end{align}
The
predictors~\eqref{eq:EDMD-dictionary-predictor}-\eqref{eq:EDMD-function-predictor}
are special cases of the approximations
in~\eqref{eq:Koopman-matrix-noninvariant}-\eqref{eq:Koopman-predictor-noninvariant},
where the orthogonal projection corresponds to the $L_2(\mu_X)$ inner
product, with empirical measure
\begin{align}\label{eq:empirical-measure}
\mu_X = \frac{1}{N} \sum_{i=1}^{N} \delta_{x_i},
\end{align}
where $\delta_{x_i}$ is the Dirac measure defined on the $i$th column
of $X$ (see e.g.,~\cite{MK-IM:18}). The quality of predictors
in~\eqref{eq:EDMD-dictionary-predictor}-\eqref{eq:EDMD-function-predictor}
depends on the quality of $\Span(\Psi)$ in terms of being close to
invariant under $\Kc$. If $\Span(\Psi)$ is invariant under $\Kc$, then
the predictors
in~\eqref{eq:EDMD-dictionary-predictor}-\eqref{eq:EDMD-function-predictor}
are exact and match
equations~\eqref{eq:Koopman-matrix-invariant}-\eqref{eq:Koopman-predictor-invariant},
resp. Determining closeness to invariance requires an
appropriate metric, which is the concept we review next.

\subsection{Consistency Index Measures The Dictionary's Quality}\label{subsec:consistency-index}
The concept of temporal forward-backward consistency to
measures how close a dictionary span is to being Koopman invariant.
Given a dictionary $\Psi$ with real-valued elements and data matrices
$X, X^+$, the \emph{consistency index}~\cite{MH-JC:23-csl}~is
\begin{align}\label{eq:consistency-index}
\ic (\Psi,X,X^+) = \lambda_{\max} (I - K_F K_B),
\end{align}
where $K_F = \psixp \psix^\dagger$ and $K_B = \psix \psixp^\dagger$
are EDMD matrices applied forward and backward in time\footnote{Note
  that this definition is equivalent but different
  from~\cite[Definition~1]{MH-JC:23-csl}. The data matrices in this
  paper are transpose of the ones in~\cite{MH-JC:23-csl}; however,
  this transposition does not affect the value of the consistency
  index.}. When the context is clear, we simply use~$\ic$.

The intuition behind the consistency index is that when $\Span(\Psi)$
is Koopman-invariant, the forward and backward EDMD matrices $K_F$ and
$K_B$ are inverse of each other. Otherwise, their product will deviate
from the identity matrix, with the consistency index providing a
measure for this deviation.  The consistency index is easy to compute
based on data and its value only depends on the vector space, not on
the choice of particular basis. The following result states a key
property of relevance to the ensuing discussion.

\begin{theorem}\longthmtitle{\cite[Theorem~1]{MH-JC:23-csl}:
$\sqrt{\ic}$ Bounds the Relative 
$L_2$-norm error for EDMD's Koopman
Predictions}\label{t:RRMSE-bound-sprad-consistency} Given
Assumption~\ref{a:full-rank} for dictionary $\Psi$, data matrices
$X,X^+$, empirical measure~\eqref{eq:empirical-measure}, consistency index~\eqref{eq:consistency-index}, and the predictor of EDMD defined in~\eqref{eq:EDMD-function-predictor}, we have
\begin{align*}
\sqrt{\ic(\Psi,X,X^+)} = \max_{f \in \Span(\Psi)} \frac{\|\Kc f -
\Pf_{\Kc f}^{\edmd}\|_{L_2(\mu_X)}}{{\| \Kc f \|_{L_2(\mu_X)}}}.
\end{align*}
The
maximum above is taken over
all functions leading to nonzero denominator (when the denominator is
zero, the numerator is also zero and the prediction is exact).
\oprocend
\end{theorem}
\smallskip

\begin{remark}\longthmtitle{Properties of the Consistency Index and
Implications for Learning} {\rm The consistency
index~\eqref{eq:consistency-index} provides a notion of
worst-case error bound for Koopman predictions on the vector
space $\Span(\Psi)$. Note that $\ic(\Psi, X, X^+) \in [0,1]$ for
all subspaces and data. Moreover, the index does not depend on
the choice of basis for the subspace: for example, if $\Phi$
provides a different basis for $\Span(\Psi)$, we have
$\ic(\Psi, X, X^+) = \ic(\Phi, X, X^+)$. Therefore, one can use
the consistency index as an effective loss function for subspace
learning\footnote{The residual error of EDMD
$\| \psixp - \Kedmd \psix \|_F$ depends on the choice of basis
and is not suitable for measuring quality of $\Span(\Psi)$. In
fact, it is easy to show~\cite[Example~1]{MH-JC:23-csl} that,
if $\Span(\Psi)$ is not invariant but contains one exact
eigenfunction, then one can find a linear transformation on
the dictionary to make the residual error \emph{arbitrarily}
close to zero.}. We refer the reader to~\cite[Proposition~1
and Lemma~1]{MH-JC:23-csl} for more information. \oprocend }
\end{remark}

\section{Motivation and Problem Statement}\label{sec:problem-statement}
Consider the discrete-time control system
\begin{align}\label{eq:control-system}
x^+ = \Tc(x,u), \quad x \in \Xc \subseteq \real^n, \; u \in \Uc
\subseteq \real^m, 
\end{align}
where $x$ and $u$ are the state and input vectors, $\Xc$ and $\Uc$
are the state and input spaces, and $\Tc: \Xc \times \Uc \to \Xc$ is the function describing the system dynamics. Note that no special
structure (e.g., control affine) is assumed on the
system~\eqref{eq:control-system}. Our goal is to provide a Koopman
operator theory description of the nonlinear control system. The
challenge for extending the concept of the Koopman operator to 
systems with inputs is that unlike 
system~\eqref{eq:dynamical-sys}, the behavior of the control
system~\eqref{eq:control-system} cannot be determined without
knowledge of the input sequence\footnote{Given an infinite input
sequence, one can determine~\cite{MK-IM-automatica:18} the system
behavior completely and define a Koopman operator for it. Moreover,
if one closes the loop by means of feedback, the system takes the
\new{form of}~\eqref{eq:dynamical-sys} and hence has a well-defined
Koopman operator, see~e.g.,~\cite{EK-JNK-SLB:21}.}.  Here, we aim to
provide a mathematical description of how to employ the Koopman
operator for control systems in both infinite- and finite-dimensional
cases, and articulate its application in data-driven modeling of
control systems.

\new{
\begin{problem}\longthmtitle{Challenges Regarding the Extension of
Koopman Theory to Control Systems}\label{prob:statement}
We aim to provide an operator-theoretic framework based on the Koopman
operator with the following properties:
\begin{enumerate}
\item \textbf{Capturing the behavior of control
    system~\eqref{eq:control-system}}: given an appropriate function
  space, the framework must encode the system behavior as evolution of
  functions similarly to~\eqref{eq:Koopman-def}
  and~\eqref{eq:multistep-Koopman-nocontrol}. This evolution must hold
  for all functions in the function space and all input sequences;
\item \textbf{Compatible with the Koopman operator for systems without
    input:} if we set the input to be constant, the framework should
  reduce to the Koopman operator for systems without
  input~\eqref{eq:Koopman-def};
\item \textbf{General finite-dimensional form:} the framework should
  provide a generalized notion of subspace invariance leading to a
  general finite-dimensional form with the following properties:
  \begin{enumerate}[label=(\roman*)]
  \item \textbf{compatibility with existing forms in the literature:}
    the general finite-dimensional form must encompass commonly used
    linear, bilinear, and linear switched control models;
  \item \textbf{compatibility with the case of systems without input
      in~\eqref{eq:Koopman-matrix-invariant}
      and~\eqref{eq:Koopman-matrix-noninvariant}:} if we set the input
    to be constant, the finite-dimensional form should reduce to the
    lifted linear model in~\eqref{eq:Koopman-matrix-invariant}
    and~\eqref{eq:Koopman-matrix-noninvariant};
  \end{enumerate}
\item \textbf{Best (optimal) approximations for finite-dimensional
    forms:} the framework should provide a notion of best
  approximation for finite-dimensional forms given a subspace of
  choice to perform the approximation;
\item \textbf{Accuracy bounds on finite-dimensional approximations:}
  the framework should provide a bound on the approximation error for
  the aforementioned models given all functions in the
  finite-dimensional subspace;
\item \textbf{Applications to data-driven modeling:} the framework
  should readily be useful for data-driven identification of control
  systems.
\end{enumerate}
\end{problem}
}

\section{Koopman Control Family and General Form for
Finite-Dimensional Models}
Here, we take the first step towards extending Koopman theory to the
control system case and providing a generalized notion of subspace
invariance. As we show below, this ultimately leads to a
finite-dimensional model that is the extension of the linear form
in~\eqref{eq:Koopman-matrix-invariant}.  We start from the observation
that, if we fix the input as a constant, we get 
\new{a system}
in the form of~\eqref{eq:dynamical-sys} which admits a well-defined
Koopman operator. Motivated by this idea, one can model the
system~\eqref{eq:control-system} by switching between constant input
systems\footnote{The idea of modeling a control system via
constant input systems has already been considered several times in
the literature~\cite{AS-DE:17,AS-AM-DE:18,SP-SK:19}.  }.  Formally,
consider the family of 
systems created by setting the input
as a constant signal
\begin{align}\label{eq:switched-family}
x^+ = \Tc_{u^*} (x) := \Tc(x, u \equiv u^*), \quad u^* \in \Uc.
\end{align}
Note that any trajectory $\{x_k\}_{k=0}^L \subset{\Xc}$ of
system~\eqref{eq:control-system} generated with input sequence
$\{u_k\}_{k=0}^{L-1} \subset{\Uc}$, can be generated by applying the
systems $\Tc_{u_k}$, $k \in \{0,\ldots,L-1\}$ subsequently
on the initial condition $x_0$. Hence, we have
\begin{subequations}
\begin{align}
x_k
&= \Tc_{u_{k-1}} \circ \cdots \circ \Tc_{u_0}
(x_0) \label{eq:mulitstep-control-family}
\\
&= \Tc_{u_{k-1}} (x_{k-1}), \quad k \in
\until{L}. \label{eq:switched-constant-input-evolution} 
\end{align}  
\end{subequations}
Noting that the members of the family $\{\Tc_{u^*}\}_{u^* \in \Uc}$
are all 
\new{without input,}
we can define Koopman operators for each of them
as in~\eqref{eq:Koopman-def}, leading to the following definition.

\begin{definition}\longthmtitle{Koopman Control Family (KCF)}\label{def:KCF}
Let $\Fc$ be a vector space (over the field $\cplx$) of
complex-valued functions with domain $\Xc$ that is closed under
composition with members of $\{\Tc_{u^*}\}_{u^* \in \Uc}$.  The
associated \textbf{Koopman control family (KCF)} is the family of
operators $\{\Kc_{u^*} : \Fc \to \Fc\}_{u^* \in \Uc}$ where, for
each $u^* \in \Uc$, $\Kc_{u^*}$ defined by
$ \Kc_{u^*} f = f \circ \Tc_{u^*}$, for all $f \in \Fc $, is the
Koopman operator corresponding to the dynamics~$\Tc_{u^*}$.
\oprocend
\end{definition}

Similarly to the multi-step
prediction~\eqref{eq:multistep-Koopman-nocontrol} under the Koopman
operator of 
\new{a system without input,}
one can
use~\eqref{eq:mulitstep-control-family} and the definition of KCF to
provide a similar identity for the 
\new{case of control systems,}
\begin{align*}
f(x_k) = \Kc_{u_0} \Kc_{u_1} \ldots \Kc_{u_{k-1}} f (x_0), \quad
\forall f \in \Fc. 
\end{align*}
Note that the identity above is \emph{exact} and general and can be
utilized for all trajectories of~\eqref{eq:control-system}, \new{thus
  answering Problem~\ref{prob:statement}(a). Moreover, if we set the
  input equal to a constant for all time to generate a system without
  input, the equation above reduces
  to~\eqref{eq:multistep-Koopman-nocontrol} and
  Definition~\ref{def:KCF} will contain a single operator as
  in~\eqref{eq:Koopman-def}.  This is in accordance with the
  requirement of Problem~\ref{prob:statement}(b).}

Even though a KCF on an appropriate function space can completely
capture the behavior of control system~\eqref{eq:control-system}, the
infinite-dimensional nature of the function space $\Fc$ can make its
implementation on digital computers impossible. To address this issue,
we need finite-dimensional representations for KCF \new{as laid out in
  Problem~\ref{prob:statement}(c)}. A simple guiding observation in
this regard is that if we have an exact finite-dimensional
representation and fix the input to be constant, \new{then the reduced
  system will have the form~\eqref{eq:dynamical-sys}}
and the model should reduce to a linear
finite-dimensional case similar to
equation~\eqref{eq:Koopman-matrix-invariant}. This leads us to the
concept of common invariant subspaces under KCF.

\begin{definition}\longthmtitle{Common Invariant Subspaces Under the
Koopman Control Family}
The space $\Lc \subseteq \Fc$ is a \textbf{common invariant
subspace} under the KCF if $\Kc_{\bar{u}} f \in \Lc$, for all
$\Kc_{\bar{u}} \in \{\Kc_{u^*}\}_{u^* \in \Uc}$ and all $f \in \Lc$.
\oprocend
\end{definition}

Finite-dimensional common invariant subspaces under the KCF
$\{\Kc_{u^*}\}_{u^* \in \Uc}$ are of utmost importance because the
action of all its members on such subspaces can be captured
\emph{exactly} by matrices. This provides a general framework for
treating control systems. Next, we show that finite-dimensional common
invariant subspaces under KCF lead to a universal form of models that
can be viewed as a generalization
of~\eqref{eq:Koopman-matrix-invariant} to the case of control
systems\footnote{Similarly to the case of 
\new{systems without input}
(cf.~Remark~\ref{r:invariance-general-form} and its preceding
discussions), we rely on a notion of subspace invariance to find a
general finite-dimensional form. In the following sections, we also
rigorously investigate approximations on non-invariant subspaces.}.

\begin{theorem}\longthmtitle{General Form on Common Invariant
Subspaces: \textbf{Input-State
Separable}}\label{t:common-invariant-separates}
The Koopman control family has a finite-dimensional (of dimension
$s$) common invariant subspace if and only if there exist functions
$\Psi:\Xc \to \cplx^s$ and $A: \Uc \to \cplx^{s \times s}$ such that
for all $(x,u) \in \Xc \times \Uc$,
\begin{align}\label{eq:composition-to-separation}
\Psi(x^+) = \Psi \circ \Tc (x,u) = \Ac(u) \Psi(x).
\end{align}
In this formulation, the common invariant subspace under the KCF is
described by $\Span(\Psi)$.
\end{theorem}
\begin{pf}
$(\Rightarrow):$ Let $\Sc \subset \Fc$, with $\dim \Sc = s$, be a
common invariant subspace of the Koopman control family
$\{\Kc_{u^*}\}_{u^* \in \Uc}$.  Let functions
$\{\psi_1, \ldots, \psi_s\}$ be a basis for $\Sc$ and define the
vector-valued function $\Psi: \Xc \to \cplx^s$ as
$\Psi(x) = [\psi_1(x),\ldots, \psi_s(x)]^T$ for all $x \in
\Xc$. Since $\Sc = \Span(\Psi)$ is invariant under the KCF, for
\emph{each} $u^* \in \Uc$, there exists a matrix
$K_{u^*} \in \cplx^{s \times s}$ (which represents the action
of operator $\Kc_{u^*}$ in KCF on subspace $\Sc$ with respect to
basis $\Psi$)
such that 
\begin{align}\label{eq:operator-by-matrix}
\Kc_{u^*} \Psi(x) = K_{u^*} \Psi (x), \; \forall x \in \Xc, 
\end{align}
where we have used~\eqref{eq:Koopman-matrix-invariant} and the
notation in Remark~\ref{r:notation-overload}.  Define then the 
matrix-valued function $\Ac: \Uc \to \cplx^{s \times s}$ as
\begin{align*}
\Ac(u^*) = K_{u^*},
\end{align*}
for each $u^* \in \Uc$.  Noting that
equation~\eqref{eq:operator-by-matrix} holds for all
$u^* \in \Uc$, one can use the definition of
$\Ac: \Uc \to \cplx^{s \times s}$ and write
\begin{align*}
\Psi \circ \Tc(x,u^*) = K_{u^*} \Psi(x) = \Ac(u^*) \Psi(x),
\; \forall x \in \Xc, \, \forall u^* \in \Uc.
\end{align*}
Noting that the equation above holds for all $u^*\in
\Uc$, one can do a change of symbol ($u^*$ to $u$), leading
to~\eqref{eq:composition-to-separation}.

$(\Leftarrow):$ Assume equation~\eqref{eq:composition-to-separation}
holds. Hence, for all $u^* \in \Uc$,
\begin{align*}
\Psi \circ \Tc(x, u \equiv u^*) = \Psi \circ \Tc_{u^*} (x) = \Ac(u
\equiv u^*) \Psi(x), \forall x \in \Xc. 
\end{align*}
Given that $\Ac (u \equiv u^*)$ is a constant matrix,
for any function $f \in \Span(\Psi)$ in the form of
$f = v_f^T \Psi$ with $v_f \in \cplx^s$, we have
\begin{align*}
\Kc_{u^*} f = f \circ \Tc_{u^*} = v_f^T \Psi \circ \Tc_{u^*} =
v_f^T \Ac(u \equiv u^*) \Psi \in \Span(\Psi). 
\end{align*}
This equality holds for all $u^* \in \Uc$; hence, $\Span(\Psi)$
is invariant under the Koopman control family
$\{\Kc_{u^*}\}_{u^* \in \Uc}$. \qed
\end{pf}

The input-state separable form~\eqref{eq:composition-to-separation}
(note the composition on the left and the matrix product on the right)
can be viewed as a generalization
of~\eqref{eq:Koopman-matrix-invariant}, which describes the exact
action of the Koopman operator on an invariant subspace.  Next,
we discuss the special structure of the input-state separable models
and their closed-form solutions.
\begin{remark}\longthmtitle{Linearity in Lifted State and
Closed-Form Solution of Input-State Separable Models} {\rm By
defining the lifting map $x \mapsto z:=\Psi(x)$, one can rewrite
the input-state separable form
in~\eqref{eq:composition-to-separation} as $z^+ = \Ac(u) z$.
This system is linear in the lifted state $z$ and has a
closed-from solution: given initial condition $x_0$ and input
sequence $\{ u_i \}_{i=0}^\infty$, one can write
\begin{align*}
z_k = \Big( \prod_{i=0}^{k-1} \Ac(u_i) \Big) z_0, \quad
\forall k \in \naturals, 
\end{align*}
where
$\prod_{i=0}^{k-1} \Ac(u_i) := \Ac(u_{k-1}) \Ac(u_{k-2}) \cdots
\Ac(u_{0})$ and $z_0 = \Psi(x_0)$. \oprocend }
\end{remark}

It is important to note that the condition in
Theorem~\ref{t:common-invariant-separates} is \emph{necessary and
sufficient (and hence cannot be relaxed); therefore the
input-state separable form is general.} In fact, as we show next, it
provides a mathematical framework encompassing common Koopman-inspired
descriptions of the control system~\eqref{eq:control-system}.
It is easy to see that the linear switched systems used
in~\cite{SP-SK:19} are a special case of the input-state separable
form where the input space $\Uc$ contains finitely many elements. We
formalize this observation in the following result that follows
directly from the definition of input-state separable form.

\begin{lemma}\longthmtitle{Linear Switched Form is a Special Case of
Input-State Separable Form}\label{l:switched-linear-separable}
For system~\eqref{eq:control-system}, let $\Uc = \{u_1,\ldots,u_l\}$
and assume the system has a lifted linear switched
representation of the form
\begin{align}\label{eq:linear-switched}
\Psi(x^+) \!=\! A_u \Psi(x), \; A_u \in \{A_{u_1}, \dots, A_{u_l}
\}\subset \real^{N_\Psi \times N_\Psi}, 
\end{align}
where $\Psi: \Xc \to \real^{N_{\Psi}}$ with $N_\Psi \in \naturals$
and $u \in \Uc$.
Then, $\Span(\Psi)$ is a finite-dimensional common invariant
subspace under the KCF associated with the system
and~\eqref{eq:linear-switched} is an input-state separable
representation.  \oprocend
\end{lemma}
\smallskip

Next, we show that the commonly used linear and bilinear Koopman-based
models are also special cases of the input-state separable form.

\begin{lemma}\longthmtitle{Linear and Bilinear Forms are Special Cases
of Input-State Separable Form}\label{l:linear-bilinear-separable}
Assume the system~\eqref{eq:control-system} has a finite-dimensional
lifted representation of the form
\begin{align}\label{eq:2nd-order-model}
\psi(x^+) = A \psi(x) + \sum_{i=1}^m  B_i \psi(x) u_i + C u,
\end{align}
where $\psi: \Xc \to \real^{N_{\psi}}$ with $N_\psi \in
\naturals$. Moreover, $A, B_i \in \real^{N_\psi \times N_\psi}$ for
$i \in \until{m}$ and $C \in \real^{N_\psi \times m}$ where $m$ is
the dimension of the input vector. 
Then $\Span(\psi) + \Span(1_\Xc)$ is a finite-dimensional
common invariant subspace under the KCF associated with the
system\footnote{Here, $1_\Xc: \Xc \to \cplx$ is the constant
function defined by $1_\Xc(x) = 1$ for all $x \in \Xc$.},
which has the input-state separable representation
\begin{align}\label{eq:linear-bilinear}
\begin{bmatrix}
\psi(x^+)
\\
1_\Xc(x^+)
\end{bmatrix}
= \begin{bmatrix}
A+\sum_{i=1}^m  u_i B_i
& Cu
\\
0
&
1
\end{bmatrix}
\begin{bmatrix}
\psi(x)
\\
1_\Xc(x)
\end{bmatrix}.
\end{align}
\end{lemma}
\begin{pf}
Using the constant function $1_\Xc$, one can rewrite the
dynamics~\eqref{eq:2nd-order-model} as~\eqref{eq:linear-bilinear},
which is in input-state separable form. Therefore, based on
Theorem~\ref{t:common-invariant-separates},
$\Span(\psi) + \Span(1_\Xc)$ is a finite-dimensional common
invariant subspace under the KCF associated with the system. \qed
\end{pf}

\begin{remark}\longthmtitle{Existence of Linear or Bilinear Forms
Implies Common Invariant Subspaces of KCF} {\rm
Lemma~\ref{l:linear-bilinear-separable} shows that if a system has
a linear or bilinear lifted form, then its associated KCF has a
common invariant subspace. However, the converse does not hold, as
corroborated by the necessary and sufficient condition in
Theorem~\ref{t:common-invariant-separates}. Therefore, for a
system to have a linear or bilinear lifted form, stronger
conditions than the existence of common invariant subspace under
KCF are required. \oprocend }
\end{remark}

Note that linear and bilinear models are special cases of the model
in~\eqref{eq:2nd-order-model}. Therefore, the input-state separable
model captures these important special cases.

\begin{example}\longthmtitle{Input-State Separable
Form}\label{ex:input-states-separable}
{\rm
Consider
\begin{align}\label{eq:ex-control-system}
x_1^+
&= a x_1 +b u
\nonumber
\\
x_2^+
&= c x_2 + d x_1^2  + e x_1 u + fu + g \sin(u) + h
\end{align}
where $x_1, x_2$ are the state variables and $u$ is the input. The
system has the input-state separable form
\begin{align}\label{eq:ex-input-state-separable-form}
\begin{bmatrix}
x_1
\\
x_2
\\
x_1^2
\\
1
\end{bmatrix}^+
\!\!= 
\begin{bmatrix}
a & 0 & 0 & b \,u 
\\
e \, u & c & d & fu + g \sin(u) + h 
\\
2ab \, u & 0 & a^2 & b^2 u^2 
\\
0 & 0 & 0 & 1
\end{bmatrix}
\begin{bmatrix}
x_1
\\
x_2
\\
x_1^2
\\
1
\end{bmatrix}.
\end{align}
Note that for any constant $u\equiv u^*$, the system turns into an
exact lifted linear form on a Koopman invariant subspace (compare
with the linear switched model in~\cite{SP-SK:19} and
Lemma~\ref{l:switched-linear-separable}). If $b = g = 0$,
\eqref{eq:ex-input-state-separable-form} turns into the following
bilinear form (cf.~Lemma~\ref{l:linear-bilinear-separable})
\begin{align*}
\begin{bmatrix}
x_1 \\ x_2 \\ x_1^2 \\ 1
\end{bmatrix}^+
=  
\begin{bmatrix}
a & 0 & 0 & 0 
\\
0 & c & d & h 
\\
0 & 0 & a^2 & 0
\\
0 & 0 & 0 & 1
\end{bmatrix}
\begin{bmatrix}x_1 \\ x_2 \\ x_1^2 \\ 1 \end{bmatrix}
+
\begin{bmatrix}
0 & 0 & 0 & 0 
\\
e & 0 & 0 & f 
\\
0 & 0 & 0 & 0
\\
0 & 0 & 0 & 0
\end{bmatrix}
\begin{bmatrix}x_1 \\ x_2 \\ x_1^2 \\ 1 \end{bmatrix} u.
\end{align*}
If in addition we have $e = 0$, the previous equation can be
written in linear form. \oprocend
}
\end{example}
\smallskip

So far, we have established the KCF modeling can completely capture
the behavior of the control
system~\eqref{eq:control-system}. Moreover, we have found the general
form of finite-dimensional models on the common invariant subspaces
of~KCF. However, given that, in general, the KCF contains uncountably
many linear operators, one needs to find tractable ways to find or
approximate finite-dimensional common invariant subspaces under the
KCF. We tackle this task in the following sections.

\section{Parameterizing the Koopman Control Family}
We provide a way to parametrize a Koopman Control Family via a single
linear operator defined on an augmented function space. This allows us
to provide an equivalent characterization for a common invariant
subspace under the KCF.

\subsection{Augmented Koopman Operator}
To parametrize the KCF, we first parametrize the family of constant
input systems in~\eqref{eq:switched-family} as the
following augmented dynamical system
\begin{align*}
\begin{bmatrix}
x
\\
u
\end{bmatrix}^+
=
\begin{bmatrix}
\Tc(x,u)
\\
u
\end{bmatrix}.
\end{align*}
For convenience, we define the following tuple notation for the system
above
\begin{align}\label{eq:augmented-system}
(x^+,u^+) =  \Taug(x,u) := (\Tc(x,u), u), 
\end{align}
for $(x,u) \in \Xc \times \Uc$.  Note that
in~\eqref{eq:augmented-system}, $u$ is a part of the state vector and
not an input. The next result shows that this augmented system
captures the behavior of all members of constant-input systems defined
in~\eqref{eq:switched-family}.

\begin{lemma}\longthmtitle{Augmented System Parametrizes the
Constant-Input Family}\label{l:augsystem-constantinputs}
For the augmented system~\eqref{eq:augmented-system}, it
holds:
\begin{enumerate}
\item the set $\Xc \times \{u^*\}$ is forward invariant
under~\eqref{eq:augmented-system} for all $u^* \in \Uc$;
\item for any $u^* \in \Uc$, let $\{x_i\}_{i=1}^\infty$ be a
trajectory of $\Tc_{u^*}$ in~\eqref{eq:switched-family} with
initial condition $x_0 \in \Xc$ and let
$\{ (x_i^{\aug}, u_i^{\aug}) \}_{i=1}^\infty$ be a trajectory of
$\Taug$ starting from $(x_0^{\aug}, u^*) \in \Xc \times \Uc$ with
$x_0^{\aug} = x_0$. Then, $x_i = x_i^{\aug}$ for all
$i \in \naturals$. \oprocend
\end{enumerate}
\end{lemma}
%

\smallskip The proof of Lemma~\ref{l:augsystem-constantinputs}
directly follows from the definition of
system~\eqref{eq:augmented-system} and is omitted for space reasons.
As a result of Lemma~\ref{l:augsystem-constantinputs}(a), if we
restrict the state space of~\eqref{eq:augmented-system} to
$\Xc \times \{u^*\}$ for any $u^* \in \Uc$, we get a well-defined
dynamics, which we denote by $\restr{\Taug}{\Xc \times
  \{u^*\}}$. Moreover, by Lemma~\ref{l:augsystem-constantinputs}(b),
$\restr{\Taug}{\Xc \times \{u^*\}}$ captures the behavior of
$\Tc_{u^*}$ for all $u^* \in \Uc$. It is in this sense we say $\Taug$
on the state space $\Xc \times \Uc$ parametrizes the family of
constant-input systems $\{\Tc_{u^*}\}_{u^* \in \Uc}$.

Since the augmented system 
\new{does not have input and is in the form~\eqref{eq:dynamical-sys},}
we can define a Koopman operator as given
in~\eqref{eq:Koopman-def}.  Appropriately defined, this operator would
encompass the KCF's information, as supported by
Lemma~\ref{l:augsystem-constantinputs}, which connects the augmented
system~\eqref{eq:augmented-system} to the constant-input
systems~\eqref{eq:switched-family}. Nonetheless, before making this
connection, we must first bridge the gap between the state-space of
constant-input systems ($\Xc$) and that of the augmented system
($\Xc \times \Uc$), and define a proper function space. To do this, we
first provide the following definition.

\begin{definition}\longthmtitle{Control-Independent Extension of
Functions in $\Fc$ to Domain $\Xc \times
\Uc$}\label{def:cont-indep-extenstion} 
Given the function $\phi \in \Fc$ where $\phi: \Xc \to \cplx$, we
define its \textbf{control-independent extension} to the domain
$\Xc \times \Uc$ as $\phi_e: \Xc \times \Uc \to \cplx$,
\begin{align*}
\phi_e(x,u) = \phi(x) 1_\Uc(u) , \quad \forall (x,u) \in \Xc
\times \Uc ,     
\end{align*}
where $1_\Uc: \Uc \to \cplx$ is defined as $1_\Uc(u) = 1$ for all
$u \in \Uc$. Similarly, for a vector-valued function
$\Phi(x) = [\phi_1(x), \ldots, \phi_n(x)]^T$, where $\phi_i \in \Fc$
for all $i \in \until{n}$, we define
$\Phi_e(x,u) = [\phi_1(x)1_\Uc(u), \ldots, \phi_n(x)1_\Uc(u)]^T$.
\oprocend
\end{definition}
\smallskip

One could equivalently define the control-independent extension as
$\phi_e(x,u) = \phi(x)$ for $(x,u) \in \Xc \times \Uc$. However, the
structure of Definition~\ref{def:cont-indep-extenstion} is consistent
with input-state separable forms, which is particularly convenient in
our forthcoming theoretical analysis.  Next, we state straightforward
but useful properties of control-independent extensions that follow
from the definition.

\begin{lemma}\longthmtitle{Control-Independent Extensions'
Properties}\label{l:control-indep-properties}
Let $\phi: \Xc \to \cplx$ and $\Phi: \Xc \to \cplx^n$, and let
$I_\Uc^{n \times n}: \Uc \to \cplx^{n \times n}$ be a constant
function returning the identity matrix, i.e.,
$I_\Uc^{n \times n}(u) = I_{n \times n}$ for all $u \in
\Uc$. Then, for all $(x,u) \in \Xc \times \Uc$,
\begin{enumerate}
\item $\phi_e(x,u) = \phi(x)$ and $\Phi_e(x,u) = \Phi(x)$;
\item $\Phi_e(x,u) = I_\Uc^{n \times n}(u) \Phi(x)$;
\item for all $f \in \Span(\Phi)$ with description
$f = v_f^T \Phi$ where $v_f \in \cplx^n$, we have
$f_e = v_f^T \Phi_e$.  \oprocend
\end{enumerate}
\end{lemma}
\smallskip

We next define a proper function space for the Koopman operator
associated with the augmented system~\eqref{eq:augmented-system}.

\begin{definition}\longthmtitle{Function Space and Koopman
Operator for $\Taug$}\label{def:augmented-operator}
Let $\Faug$
be a linear space (on the field $\cplx$) of complex-valued
functions with domain $\Xc \times \Uc$ such that
\begin{enumerate}
\item is closed under composition with $\Taug$;
\item contains $ f \circ \Tc$ for all $f \in \Fc$;
\item contains the control-independent extension $f_e$ for all
$f \in \Fc$; 
\item for all $u^* \in \Uc$, $\restr{\Faug}{u=u^*} = \Fc$.
\end{enumerate}
Then, the \textbf{augmented Koopman operator},
$\Kaug: \Faug \to \Faug$ defined as
\begin{align}\label{eq:augmented-Koopman}
\Kaug g = g \circ \Taug, \quad \forall g \in \Faug,
\end{align}
encodes the behavior of the augmented
system~\eqref{eq:augmented-system}.  \oprocend
\end{definition}

Note that, as long as we allow the function spaces $\Fc$ and $\Faug$
to be infinite-dimensional, the conditions in
Definition~\ref{def:augmented-operator} are easy to satisfy.  The
choice of $\Fc$ and $\Faug$ depends on $\Tc$ and $\Taug$ and the
choices are not unique. Here, we provide a few generic examples for
$\Fc$ and $\Faug$:
\begin{itemize}
\item one of the simplest examples is to choose $\Fc$ and $\Faug$ as
the spaces of bounded complex-valued functions on domains $\Xc$ and
$\Xc \times \Uc$, resp.;
\item if $\Tc$ in~\eqref{eq:control-system} is continuous in both
variables (e.g., with respect to usual metrics inherited from
$\real^n$ and $\real^m$), then $\Fc$ and $\Faug$ can be chosen as
the spaces of continuous complex-valued functions on domains $\Xc$
and $\Xc \times \Uc$, resp.;
\item if $\Tc$ is polynomial in both variables, then $\Fc$ and $\Faug$
can be chosen as the spaces of all polynomials with complex
coefficients on $\Xc$ and $\Xc \times \Uc$, resp.
\end{itemize}
It is also important to note that, at this point of the exposition,
there are no requirements on these function spaces having inner
products, norms, or metrics.

\subsection{Augmented Koopman Operator Parametrizes the Koopman
Control Family}

Here, we investigate the connection between the augmented operator and
the KCF, and the implications for the search of common invariant
subspaces for the KCF.  The next result shows how $\Kaug$
parameterizes the KCF $\{\Kc_{u^*}\}_{u^* \in \Uc}$.

\begin{lemma}
\longthmtitle{$\Kaug$ Parametrizes the
KCF}\label{l:augmentedKoopman-captures-individuals}
Let $f \in \Faug$. Then for all $u^* \in \Uc$ we have
$\restr{(\Kaug f)}{u=u^*} = \Kc_{u^*} (\restr{f}{u=u^*})$.
\oprocend
\end{lemma}
\begin{pf}
  By definition, $\Kaug f(x,u) = f(\Tc(x,u),u)$.~Hence,
  \begin{align*}
    \restr{(\Kaug f(x,u))}{u=u^*}
    & =
      f(\Tc(x,u^*), u^*)  
      = f(\Tc_{u^*}(x), u^*)
    \\
    &
      = \restr{f}{u=u^*} (\Tc_{u^*}(x)) = \Kc_{u^*} [\restr{f}{u=u^*}] (x),
  \end{align*}
  for all $u^* \in \Uc$ and all $x \in \Xc$. \qed
\end{pf}

Lemma~\ref{l:augmentedKoopman-captures-individuals}
establishes the important fact that the action of $\Kaug$ on $\Faug$
completely captures the effect of $\Kc_{u^*}$ on
$\restr{\Faug}{u=u^*} = \Fc$
(cf.~Definition~\ref{def:augmented-operator}) for all $u^* \in
\Uc$. This shows that $\Kaug$ can be viewed as a parametrization of
the KCF, i.e., by knowing the effect of $\Kaug$ on $\Faug$, one can
calculate the effect of all (potentially uncountably many) members of
the KCF. The next result shows how the augmented Koopman operator can
capture relevant information regarding the evolution of functions in
$\Fc$ under the trajectories of the control
system~\eqref{eq:control-system}.

\begin{lemma}\longthmtitle{Augmented Koopman Operator Predicts the
Functions Evolutions on  System
Trajectories}\label{l:Kaug-captures-function-evolutions-in-F} 
Let $f \in \Fc$ and denote by $f \circ \Tc \in \Faug$ the function
created by pushing the values of $f$ one time-step forward through
the trajectories of $\Tc$.  Let $f_e$ be the control-independent
extension of $f$ to $\Xc \times \Uc$.  Then,
$f \circ \Tc = \Kaug f_e$.
\oprocend
\end{lemma}
\begin{pf}
  For all $(x,u) \in \Xc \times \Uc$, one can write
  \begin{align*}
    \Kaug f_e (x,u)
    &= f_e(\Taug(x,u)) = f_e (\Tc(x,u), u) 
    \\
      & = f(\Tc(x,u)) = f \circ \Tc (x,u),
  \end{align*}
  where we have used~\eqref{eq:augmented-system} in the second
  equality and Lemma~\ref{l:control-indep-properties}(a) in the third
  equality. \qed
\end{pf}

 Lemma~\ref{l:Kaug-captures-function-evolutions-in-F}
provides a crucial tool to analyze the behavior of functions in $\Fc$
on the trajectories of the control system~\eqref{eq:control-system}
(note the similarity of the composition $f \circ \Tc$ with the
definition of the Koopman operator~\eqref{eq:Koopman-def} for
systems \new{without input}). In this result, observe that even though $\Kaug$
is the Koopman operator associated with~\eqref{eq:augmented-system},
its action on control-independent function extensions leads to the
prediction of the function values on trajectories of the actual
control system~\eqref{eq:control-system}.

The next result provides a link between the invariant subspaces of
$\Kaug$ and common invariant subspaces of the KCF.

\begin{proposition}\longthmtitle{Invariant Subspaces of $\Kaug$
Characterize Common Invariant Subspaces for the
KCF}\label{p:family-augmented-invariance} 
Let $\Sc \subseteq \Faug$ be an invariant subspace under
$\Kaug$. Then,
\begin{enumerate}
\item for all $u^* \in \Uc$, $\restr{\Sc}{u=u^*}$ is an invariant
subspace of $\Kc_{u^*}$;
\item if $\restr{\Sc}{u=u_1} = \restr{\Sc}{u=u_2}$ for all
$u_1, u_2 \in \Uc$, then $\restr{\Sc}{u=u^*}$ (for any
$u^* \in \Uc$) is a common invariant subspace under the Koopman
control family $\{\Kc_{u^*}\}_{u^* \in \Uc}$.
\end{enumerate}
\end{proposition}
\begin{pf}
(a) First note that $\restr{\Sc}{u=u^*}$ is a vector space for all
$u^* \in \Uc$. Given any $u^* \in \Uc$, consider an arbitrary
function $g \in \restr{\Sc}{u=u^*}$. By definition of
$\restr{\Sc}{u=u^*}$, there exists a function $\tilde{g} \in \Sc$
such that $g = \restr{\tilde{g}}{u = u^*}$ (note that $\tilde{g}$
might not be unique). By
Lemma~\ref{l:augmentedKoopman-captures-individuals}, one can write
\begin{align*}
\Kc_{u^*} g = \Kc_{u^*} (\restr{\tilde{g}}{u=u^*})= \restr{(\Kaug
\tilde{g})}{u=u^*} \in \restr{\Sc}{u=u^*},
\end{align*}
where we have used the fact that $\Kaug \tilde{g} \in \Sc$ because
$\Sc$ is invariant under $\Kaug$. Therefore, $\restr{\Sc}{u=u^*}$ is
an invariant subspace of~$ \Kc_{u^* }$.

(b) This is a direct consequence of part~(a) and the definition of
common invariant subspace for the KCF. \qed
\end{pf}

Proposition~\ref{p:family-augmented-invariance} provides a tool for
the identification of common invariant subspaces under KCF based on
the invariant subspaces of the augmented Koopman operator. However,
checking the condition in
Proposition~\ref{p:family-augmented-invariance}(b) requires one to
compare different vector spaces, which can be cumbersome. In the
following section, we provide more direct conditions that can be
checked easily and lead to input-state separable models, as laid out
in Theorem~\ref{t:common-invariant-separates}.

\section{Input-State Separable Forms via the Augmented Koopman
Operator} 
Here, we aim to build on
Proposition~\ref{p:family-augmented-invariance} and
Theorem~\ref{t:common-invariant-separates} to provide more specific
practical criteria to identify common invariant subspaces of the KCF
and derive input-state separable models. Based on
Theorem~\ref{t:common-invariant-separates} we know that on a common
invariant subspace, function composition with $\Tc$ leads to functions
that can be written as a linear combination of separable functions in
$x$ and $u$. For convenience, we provide the following definition.

\begin{definition}\longthmtitle{Input-State Separable Functions and
Their Linear Combinations}\label{def:separable}
A function $f \in \Faug$ is \textbf{input-state separable} if there
exist $g: \Uc \to \cplx$ and $h: \Xc \to \cplx$ such that
$f(x,u) = g(u) h(x)$ for all $x \in \Xc$ and $u \in \Uc$.  A
function $J$ is an \textbf{input-state separable combination} (or
\textbf{separable combination} for short) if it can be written as a
\emph{finite} linear combination of input-state separable functions.
\oprocend
\end{definition}

Next, we show a property of the bases for spaces of separable
combinations.

\begin{proposition}\longthmtitle{Spaces of Separable Combinations Have
Separable Bases}\label{p:separable-basis}
Let $\Sc \subseteq \Faug$ be a finite-dimensional (of dimension
$s \in \naturals$) subspace comprised of input-state separable
combinations. Then, for any arbitrary basis
$\{\phi_1,\ldots, \phi_s\}$ of $\Sc$, the vector-valued function
$\Phi(x,u) = [\phi_1(x,u), \ldots, \phi_s(x,u)]^T$ can be decomposed
as the product of two functions as follows
\begin{align}\label{eq:basis-product-decomposition}
\Phi(x,u) = G(u) H(x),  \; \forall (x,u)\in \Xc \times \Uc.
\end{align}
where $G: \Uc \to \cplx^{s \times l}$ and $H: \Xc \to \cplx^{l}$ for
some $l \in \naturals$. \oprocend
\end{proposition}
\begin{pf}
	By hypothesis, for each $i \in \until{s}$, there exists $n_i$ such
	that
	\begin{align}\label{eq:phi_decomposition}
		\phi_i (x,u) = \sum_{j_i = 1}^{n_i} p_{j_i}^i (u) q_{j_i}^i (x),
		\; \forall (x,u)\in \Xc \times \Uc. 
	\end{align}
	for some functions $p_{j_i}^i: \Uc \to \cplx$,
	$ q_{j_i}^i: \Xc \to \cplx$.  Now, consider the space
	\begin{align*}
		Q = \Span \{ q_{j_i}^i: \Xc \to \cplx \; | \; i \in \until{s}, j_i \in \until{n_i} \}.
	\end{align*}
	By construction, $Q$ is finite dimensional, with
	$l:= \dim Q \leq \sum_{i=1}^s n_i$. Let $\{h_1, \ldots, h_l\}$ be a
	basis for $Q$ and construct the vector-valued function
	$H:\Xc \to \cplx^l$ as
	\begin{align*}
		H(\cdot) = [h_1(\cdot), \ldots, h_l(\cdot)]^T.
	\end{align*}
	By construction of $Q$, all functions $q_{j_i}^i$ can be written as
	linear combinations of $\{h_1, \ldots, h_l\}$. Hence, there exist
	vectors $v_{j_i}^i \in \cplx^l$ such that
	\begin{align}\label{eq:q-based-on-H}
		q_{j_i}^i (\cdot)= (v_{j_i}^i)^T H(\cdot), \; \forall i \in \until{s}, \,  j_i \in \until{n_i}.
	\end{align}
	Now, based on~\eqref{eq:phi_decomposition}-\eqref{eq:q-based-on-H},
	for all $i \in \until{s}$, one can write
	\begin{align}\label{eq:phi-new-decomposition}
		\phi_i (x,u) = \sum_{j_i = 1}^{n_i} p_{j_i}^i (u) (v_{j_i}^i)^T
		H(x), \; \forall (x,u) \in \Xc \times \Uc. 
	\end{align}
	Defining the function $G: \Uc \to \cplx^{s \times l}$ as
	\begin{align*}
		G(u) = \begin{bmatrix}
			\sum_{j_1 = 1}^{n_1} p_{j_1}^1 (u) (v_{j_1}^1)^T
			\\
			\vdots
			\\
			\sum_{j_s = 1}^{n_s} p_{j_s}^s (u) (v_{j_s}^s)^T
		\end{bmatrix}, \; \forall u \in \Uc,
	\end{align*}
	it follows from~\eqref{eq:phi-new-decomposition} that
	$\Phi(x,u) = G(u) H(x)$. \qed
\end{pf}

With this result in place, we can show how to obtain a common
invariant subspace of the KCF using the invariant subspaces of the
augmented Koopman operator.

\begin{theorem}\longthmtitle{Rank Condition for Identification of
Common Invariant Subspaces of KCF via
$\Kaug$}\label{t:rank-common-invariant}
Let $\Sc \subseteq \Faug$ be a finite-dimensional (of dimension
$s \in \naturals$) subspace comprised of input-state separable
combinations that is invariant under~$\Kaug$ and let
$\Phi(x,u) = G(u) H(x)$ be a decomposition of a basis for $\Sc$,
where $G: \Uc \to \cplx^{s \times l}$ and $H: \Xc \to \cplx^{l}$.
If $G(u)$ has full column rank for all $u \in \Uc$, then the space
$\Hc = \Span(H)$ is a common invariant subspace for the~KCF.
\end{theorem}
\begin{pf}
Since $\Sc$ is a finite-dimensional invariant subspace under
$\Kaug$, given the basis $\Phi$, one can represent the action of
$\Kaug$ on $\Sc$ by a matrix $A \in \cplx^{s \times s}$ as
\begin{align}\label{eq:Kaug-basis-matrix}
\Kaug \Phi = A \Phi,
\end{align}
where we have used the compact notation in
Remark~\ref{r:notation-overload}.  Using this and the decomposition
$\Phi(x,u) = G(u) H(x)$, 
\begin{align*}
\Kaug \big(G(\cdot) H(\cdot)  \big) = A \, G(\cdot) H(\cdot).
\end{align*}
With this compact description,  to invoke
Proposition~\ref{p:family-augmented-invariance}(b), we need to show
that
$\restr{\big[ \Span( G(u) H(\cdot) )\big]}{u=u^*} = \Span(G(u^*) H)$
is the same for all $u^* \in \Uc$.
We show this by establishing
\begin{align}\label{eq:subspace-restriction-equality}
\Span(G(u^*) H )= \Span(H) ,  \; \forall u^* \in \Uc. 
\end{align}
To show the inclusion from left to right, consider
$g: \Xc \to \cplx$ with $g \in \Span(G(u^*) H)$.  Hence, there is a
vector $v_g \in \cplx^s$ such that
$g(\cdot) = v_g^T G(u^*) H(\cdot)$. Defining
$w_g = G(u^*)^T v_g \in \cplx^l$, one can write
$g(\cdot) = w_g^T H(\cdot) \in \Span(H)$, proving
\begin{align}\label{eq:GH-sub-H}
\Span(G(u^*) H) \subseteq \Span(H), \; \forall u^* \in \Uc. 
\end{align}
To prove the inclusion from right to left, consider
$p: \Xc \to \cplx$ with $p \in \Span(H)$. Hence, there is a vector
$v_p \in \cplx^l$ such that $p(\cdot) = v_p^T H(\cdot)$. For a given
$u^*$, we need to show that there exists a vector $w_p \in \cplx^s$
such that $p(\cdot) = w_p^T G(u^*) H(\cdot)$. In other words, we
have to show the following linear equation holds for some
$w_p \in \cplx^s$
\begin{align}\label{eq:linear-equation-wp}
G(u^*)^T w_p = v_p.
\end{align}
Given that $G(u^*)$ has full column
rank,~\eqref{eq:linear-equation-wp} always have at least one solution,
which might not be unique. Therefore,
$p(\cdot) = w_p^T G(u^*) H(\cdot) \in \Span(G(u^*) H)$ and
consequently
\begin{align}\label{eq:H-sub-GH}
\Span(H) \subseteq \Span(G(u^*) H), \; \forall u^* \in \Uc.
\end{align}
Combining~\eqref{eq:GH-sub-H} and~\eqref{eq:H-sub-GH} yields the
subspace equality~\eqref{eq:subspace-restriction-equality}. By
Proposition~\ref{p:family-augmented-invariance}(b), we conclude that
$\Span(H)= \Hc$ is a common invariant subspace for the KCF. \qed
\end{pf}

Theorem~\ref{t:rank-common-invariant} provides an algebraic rank
condition that is far easier to check than the condition in
Proposition~\ref{p:family-augmented-invariance}.

\begin{remark}\longthmtitle{A Note on Rank Condition in
Theorem~\ref{t:rank-common-invariant}}
{\rm 
In Theorem~\ref{t:rank-common-invariant}, if the matrix $G(u)$ is
column-rank deficient only for some $u \in \Uc$, one might be able
to use the result with a slight relaxation. In particular, define
\begin{align*}
\tilde{\Uc} :=\setdef{u \in \Uc}{G(u)\: \text{has full column rank}}.
\end{align*}
If the control system~\eqref{eq:control-system} exhibits favorable control
properties, e.g., controllability, reachability, or stabilizability,
etc., on $\tilde{\Uc}$, then one can restrict the input space to
$\tilde{\Uc}$ and apply Theorem~\ref{t:rank-common-invariant}. A
notable example of this restriction is the case of switched linear
modeling, see e.g.,~\cite{SP-SK:19}, that only requires $\tilde{\Uc}$
to contain finitely many predetermined inputs.  \oprocend
}
\end{remark}

\begin{example}\longthmtitle{Revisiting
Example~\ref{ex:input-states-separable} -- Invariant Subspace for
$\Kaug$}\label{ex:Kaug-invariant-decomposition}
{\rm
For the system~\eqref{eq:ex-control-system}, one can derive a
lifted linear form on an invariant subspace of $\Kaug$ as
\begin{align*}
\begin{bmatrix}
x_1 \\ x_2 \\ x_1^2 \\ 1 \\ x_1 u \\ u \\ u^2 \\ \sin(u)
\end{bmatrix}^+
= 
\begin{bmatrix} %
a & 0 & 0 & 0 & 0 & b &  0 & 0 %
\\
0 & c & d & h & e & f &  0 & g  %
\\
0 & 0 & a^2 & 0 & 2 ab & 0 &  b^2 & 0  %
\\
0 & 0 & 0 & 1 & 0 & 0 &  0 & 0  %
\\
0 & 0 & 0 & 0 & a & 0 &  b & 0  %
\\
0 & 0 & 0 & 0 & 0 & 1 &  0 & 0  %
\\
0 & 0 & 0 & 0 & 0 & 0 &  1 & 0  %
\\
0 & 0 & 0 & 0 & 0 & 0 &  0 & 1  %
\end{bmatrix}
\begin{bmatrix}
x_1 \\ x_2 \\ x_1^2 \\ 1 \\ x_1 u \\ u \\ u^2 \\ \sin(u)
\end{bmatrix}.
\end{align*}
Note that the evolution is based on the augmented
system~\eqref{eq:augmented-system}, which does not evolve $u$. One
can decompose the basis
$\Phi(x,u) = [x_1, x_2, x_1^2, 1, x_1 u, u, u^2, \sin(u)]^T$ as
$\Phi(x,u) = G(u) H(x)$, with
$G(u) = [I_{4 \times 4}, \tilde{G}(u)^T]^T$ where
\begin{align*}
H(x) =
\begin{bmatrix}
x_1
\\
x_2
\\
x_1^2
\\
1
\end{bmatrix}
\; \text{and} \; 
\tilde{G}(u) = 
\begin{bmatrix}
u & 0 & 0 & 0 
\\
0 & 0 & 0 & u 
\\
0 & 0 & 0 & u^2
\\
0 & 0 & 0 & \sin(u)
\end{bmatrix}.
\end{align*}
In this decomposition, the rank condition in
Theorem~\ref{t:rank-common-invariant} holds (note the presence of
$I_{4 \times 4}$ in $G(u)$). Hence, $\Span(H)$ is a common invariant
subspace for the KCF, which is in agreement with
Example~\ref{ex:input-states-separable}. \oprocend
}
\end{example}

According to Theorem~\ref{t:common-invariant-separates}, a common
invariant subspace for the KCF comes with an associated input-state
separable model for the control system~\eqref{eq:control-system}. The
next result specifies how to obtain it under the conditions of
Theorem~\ref{t:rank-common-invariant}.

\begin{proposition}\longthmtitle{Deriving Input-State Separable Models
using Invariant Subspaces of
$\Kaug$}\label{p:separable-form-invariant-augmented} 
Let $\Sc \subseteq \Faug$ be a finite-dimensional (of dimension
$s \in \naturals$) subspace comprised of input-state separable
combinations that is invariant under~$\Kaug$ and
$\Phi: \Xc \times \Uc \to \cplx^s$ a basis of $\Sc$. Let
$A \in \cplx^{s \times s}$ be\footnote{The existence of this matrix
is a direct consequence of the fact that $\Sc$ is invariant under
$\Kaug$.} such that $ \Kaug \Phi = A \Phi $.  Let
$\Phi(x,u) = G(u) H(x)$ be a decomposition of a basis for $\Sc$,
where $G: \Uc \to \cplx^{s \times l}$ and $H: \Xc \to \cplx^{l}$.
If $G(u)$ has full column rank for all $u \in \Uc$, then the
matrix-valued map $\Ac: \Uc \to \cplx^{l \times l}$ given by
\begin{align*}
\Ac(u) = G(u)^\dagger A G(u) =\big( G(u)^H G(u) \big) ^{-1} G(u)^H A G(u),
\end{align*}
turns the common-invariant subspace $\Hc = \Span(H)$ for the KCF
into the input-state separable form of
Theorem~\ref{t:common-invariant-separates}, i.e., for all
$(x,u) \in \Xc \times \Uc$
\begin{align*}
H(x^+) = H \circ \Tc (x,u) = \Ac(u) H(x).
\end{align*}
\end{proposition}
\begin{pf}
Using the definition of $\Taug$,
cf. equation~\eqref{eq:augmented-system}, one can write
$\Kaug \Phi (x,u)= \Phi \circ \Taug(x,u) = \Phi(\Tc(x,u),u) = A
\Phi(x,u)$ for all $(x,u) \in \Xc \times \Uc$. Now, by using
$\Phi(x,u) = G(u) H(x)$,
\begin{align*}
G(u) H(\Tc(x,u)) = A G(u) H(x), \quad \forall (x,u) \in \Xc \times \Uc.
\end{align*}
Keeping in mind that $G(u)$ has full column rank, one can multiply
both sides from the left by
$G(u)^\dagger = \big( G(u)^H G(u) \big) ^{-1} G(u)^H$, use
$G(u)^\dagger G(u) = I$, and reorder the terms to write
\begin{align*}
H \circ \Tc(x,u) = H(\Tc(x,u)) = G(u)^\dagger A G(u) H(x) = \Ac(u) H(x), 
\end{align*}
for all $(x,u) \in \Xc \times \Uc$. \qed
\end{pf}

Theorem~\ref{t:rank-common-invariant} and
Proposition~\ref{p:separable-form-invariant-augmented} provide us with
a way to leverage the augmented Koopman operator $\Kaug$ to identify
common invariant subspaces for the KCF and derive input-state
separable models for the control system.

\section{Input-State Separable Forms on Normal Spaces}
In this section, we turn our attention to a special case of subspaces
that are of practical significance. This focus is motivated
by examining the rank condition on~$G(u)$ presented in
Theorem~\ref{t:rank-common-invariant}, and observing that the
matrix-valued function $G$ constitutes an element of the basis
description for subspace~$\Sc$.  It becomes then clear that this
condition specifies a structural characteristic of the subspace $\Sc$
and its basis, which is independent of the operator $\Kaug$.
Therefore, here we study a specific class of subspaces and their bases
that always satisfy the rank condition in
Theorem~\ref{t:rank-common-invariant}.

\begin{definition}\longthmtitle{Vector-Valued Function of Separable
Combinations in Normal Form and Normal
Spaces}\label{def:normal-form-dictionary}
Let $\Phi: \Xc \times \Uc \to \cplx^s$ be a vector-valued function
of separable combinations. Moreover, let the set of elements of
$\Phi$ be linearly independent. $\Phi$ is in \textbf{normal form} if
it has a decomposition as one of the following:
\begin{subequations}\label{eq:normal-form-dictionary}
\begin{align}
&\Phi(x,u) =
\begin{bmatrix}
I_{\Uc}^{l \times l}(u)
\\
\tilde{G}(u)
\end{bmatrix}
H(x),
&&\text{ $s> l$}, 
\label{eq:normal-form-s-gt-l}
\\
&\Phi(x,u) =  I_{\Uc}^{l \times l}(u)  H(x),
&&
\text{ $s = l$},
\end{align}
\end{subequations}
where $H: \Xc \to \cplx^l$ and
$\tilde{G}: \Uc \to \cplx^{(s-l)\times l}$ for some $l \leq s$ and
the elements of $H$ are linearly independent.  Moreover,
$I_{\Uc}^{l \times l} : \Uc \to \cplx^{l \times l}$ is the constant
identity function, $I_{\Uc}^{l \times l}(u) \equiv I$.
A finite-dimensional space of separable combinations is
\textbf{normal} if it has a basis that can be written as a
vector-valued function of normal form.  \oprocend
\end{definition}

From Definition~\ref{def:normal-form-dictionary}, it is clear that a
basis in normal form satisfies the rank condition in
Theorem~\ref{t:rank-common-invariant}. The next result shows a useful
property of normal spaces.

\begin{proposition}\longthmtitle{Normal Spaces Capture Control-Independent
Functions}\label{p:normal-space-u-independent}
Let $\Sc \subset \Faug$ be a finite-dimensional space of input-state
separable combinations and let $\Phi(x,u)= G(u) H(x)$ be a
decomposition of a basis for $\Sc$, where
$G: \Uc \to \cplx^{s \times l}$ and $H: \Xc \to \cplx^{l}$ for some
$l \leq s$ (here, $s \in \naturals$ is the dimension of
$\Sc$). Then, $\Sc$ is normal if and only if $h_e \in \Sc$ (cf.~Definition~\ref{def:cont-indep-extenstion}) for all
$h \in \Span(H)$.
\end{proposition}
\begin{pf}
$(\Rightarrow):$ Since $\Sc$ is normal, it has a basis with one of
the following forms:
\begin{align}\label{eq:normal-forms}
&\hat{\Phi}(x,u) =
\begin{bmatrix}
I_{\Uc}^{l \times l} (u)
\\
\tilde{G}(u)
\end{bmatrix}
\tilde{H}(x),
&&\text{if $s> l$},
\nonumber
\\
&\hat{\Phi}(x,u) = I_{\Uc}^{l \times l} (u) \tilde{H}(x),
&&\text{if
$s = l$},
\end{align}
where $\Span(\tilde{H}) = \Span(H)$ (this is a direct consequence of
the fact that $\Phi$ and $\hat{\Phi}$ are bases for the same
subspace). Then, for every $h \in \Span(H)$, there exists a vector
$w_h \in \cplx^l$, such that $h(\cdot) = w_h^T
\tilde{H}(\cdot)$. Based on~\eqref{eq:normal-forms}, one can write
\begin{align*}
&h_e(x,u)  = [w_h^T, 0_{1 \times (s-l)}] \hat{\Phi}(x,u) \in \Sc,
&&\text{if $s> l$}, 
\\
&h_e(x,u) = w_h^T \hat{\Phi}(x,u) \in \Sc, &&\text{if $s = l$}.
\end{align*}
This concludes the proof of this part.

$(\Leftarrow):$ Let $H(\cdot) = [h_1(\cdot), \ldots,
h_l(\cdot)]^T$. By hypothesis, we have $h_i(x) 1_\Uc(u) \in \Sc$ for
all $i \in \until{l}$. As a result, there exist vectors
$\{w_1, \ldots, w_l\} \subset \cplx^s$ such that for all
$(x,u) \in \Xc \times \Uc$,
\begin{align}\label{eq:u-independent-representation}
h_i(x) 1_\Uc(u) = w_i^T G(u) H(x), \; \forall i \in \until{l}.
\end{align}
Let $W = [w_1, \ldots, w_l]^T \in \cplx^{l \times s}$ and consider
two cases:

Case~(i): Suppose $s = l$. Define the vector-valued function
$\tilde{\Phi}(\cdot,\cdot) = W \Phi(\cdot,\cdot)$. This function can
be written as
\begin{align*}
\tilde{\Phi}(x,u) = W G(u) H(x) = I_{\Uc}^{l \times l}(u) H(x), \quad
\forall (x,u) \in \Xc \times \Uc, 
\end{align*}
where we have used~\eqref{eq:u-independent-representation}. Therefore
$\tilde{\Phi}$ is a normal-form basis and hence $\Sc$ is normal.

Case~(ii) Suppose $s > l$ and decompose
$G(u) = [G_1^T(u), G_2^T(u)]^T$, where
$G_1: \Uc \to \cplx^{l \times l}$ and
$G_2: \Uc \to \cplx^{(s-l) \times l}$. Define the vector-valued
function $\hat{\Phi}: \Xc \times \Uc \to \cplx^{s \times l}$,
\begin{align*}
\hat{\Phi}(x,u) =
\begin{bmatrix}
W \\ B
\end{bmatrix}
\Phi(x,u)
& =
\begin{bmatrix}
W  \\ B
\end{bmatrix}
\begin{bmatrix}
G_1(u) \\ G_2(u)
\end{bmatrix}
H(x) 
\\
&
= \begin{bmatrix}
I_{\Uc}^{l \times l} (u) \\ G_2(u)
\end{bmatrix}
H(x),
\end{align*}
where
$B = [0_{(s-l) \times l}, I_{(s-l) \times (s-l)}] \in \cplx^{(s-l)
\times s}$, and in the third equality we have
used~\eqref{eq:u-independent-representation}. Therefore,
$\hat{\Phi}$ is a normal-form basis and hence $\Sc$ is normal. \qed
\end{pf}

Proposition~\ref{p:normal-space-u-independent} reveals a useful
property of normal spaces that allows us to directly predict the
evolution of functions in $\Fc$ by applying $\Kaug$ on
control-independent extensions through
Lemma~\ref{l:Kaug-captures-function-evolutions-in-F}, as
detailed~next.

\begin{theorem}\longthmtitle{Identification of Common Invariant
Subspaces of the KCF and Input-State Separable Forms on
Normal Spaces}\label{t:normal-space-common-invariant}
Let $\Sc \subset \Faug$ be a finite-dimensional normal space of
input-state separable combinations that is invariant under
$\Kaug$. Let $\Phi(x,u) = G(u) H(x)$ be a decomposition of a basis for
$\Sc$ where $G: \Uc \to \cplx^{s \times l}$ and
$H: \Xc \to \cplx^{l}$ for some $l \leq s$ (here, $s \in \naturals$
is the dimension of $\Sc$). Then,
\begin{enumerate}
\item $\Span(H) \subset \Fc$ is a common invariant subspace under
the Koopman control family $\{\Kc_{u^*}\}$;
\item for all $h \in \Span(H)$ and for all
$(x,u) \in \Xc \times \Uc$, it holds that
$h(x^+) = h\circ \Tc (x,u) = \Kaug h_e(x,u)$;
\item without loss of generally, assume $\Phi$ is in normal form,
i.e., $G(u) = I_{\Uc}^{l \times l}(u)$ if $l = s$ or
$G(u) = [I_{\Uc}^{l \times l}(u)^T, \tilde{G}(u)^T]^T$ if $s >l$.
Moreover, let $A \in \cplx^{s \times s}$ be a matrix such that
$\Kaug \Phi = A \Phi$ (note that $A$ exists because $\Sc$ is
invariant under $\Kaug$). If $s > l$, consider the
block-decomposition of $A$,
\begin{align*}
A = \begin{bmatrix} A_{11} & A_{12} \\ A_{21} & A_{22} \end{bmatrix}, 
\end{align*}
where $A_{11} \in \cplx^{l \times l}$,
$A_{12} \in \cplx^{l \times (s-l)}$,
$A_{21} \in \cplx^{(s-l) \times l}$, and
$A_{22} \in \cplx^{(s-1) \times (s-l)}$. Then, the associated
input-state separable dynamics can be written as
\begin{align}\label{eq:normal-separable-form}
H(x^+) = H \circ \Tc(x,u) = \Ac(u) H(x),
\end{align}
where, for each $u \in \Uc$,
\begin{align*}
&\Ac(u) = A_{11} + A_{12} \tilde{G}(u), &&\text{if $s > l$}. 
\\
&\Ac(u) = A \, I_{\Uc}^{l \times l}(u) = A, &&\text{if $s = l$}. 
\end{align*}
\end{enumerate}
\end{theorem}
\begin{pf}
(a) Since $\Sc$ is normal, one can do a linear transformation of the
basis $\Phi(x,u) = G(u) H(x)$ to put it in normal form. Hence, there
is a nonsingular square matrix $E$, such that
$E\Phi(x,u) = EG(u) H(x)$ is in normal form. Therefore, by
Definition~\ref{def:normal-form-dictionary}, $EG(u)$ has full column
rank for all $u \in \Uc$. Since $E$ is nonsingular, we deduce that
$G(u)$ has full column rank for all $u \in \Uc$. As a result, we can
invoke Theorem~\ref{t:rank-common-invariant} to deduce that
$\Span(H) \subset \Fc$ is a common invariant subspace under the KCF.

(b) This part is the direct consequence of
Lemma~\ref{l:Kaug-captures-function-evolutions-in-F}.

(c) Using the definition of $\Taug$,
cf. equation~\eqref{eq:augmented-system}, one can write
$ \Phi \circ \Taug(x,u) = \Phi(\Tc(x,u),u) = A \Phi(x,u)$ for all
$(x,u) \in \Xc \times \Uc$. Now, using $\Phi(x,u) = G(u) H(x)$,
\begin{align}\label{eq:normal-invariant-separation}
G(u) H(x^+) = G(u) H(\Tc(x,u)) = A G(u) H(x).
\end{align}
The case $s = l$ is trivial since $G(u)$ is an identity map. For the
case $s > l$, the proof directly follows by multiplying both sides
of~\eqref{eq:normal-invariant-separation} from the left by the
matrix $[I_{l \times l}, 0_{l \times (s-l)}]$ and using the
decompositions of $G(u)$ and~$A$. \qed
\end{pf}

\begin{example}\longthmtitle{Examples~\ref{ex:input-states-separable}
and~\ref{ex:Kaug-invariant-decomposition} Revisited}
{\rm
The basis decomposition in
Example~\ref{ex:Kaug-invariant-decomposition} is in normal form. One
can readily use the formula in
Theorem~\ref{t:normal-space-common-invariant}(c) with this
decomposition to calculate the input-state separable
form~\eqref{eq:ex-input-state-separable-form}.  \oprocend
}
\end{example}

Theorem~\ref{t:normal-space-common-invariant} has significant
practical implications: not only it connects the invariant subspaces
of $\Kaug$ to common invariant subspaces of the KCF, but more
importantly, unlike
Proposition~\ref{p:separable-form-invariant-augmented}, it provides a
direct way of predicting the evolution of observables in $\Fc$ under
the control system based on the application of $\Kaug$ on
control-independent extensions. This direct computation does not
require taking a pseudo-inverse (cf.
Proposition~\ref{p:separable-form-invariant-augmented}) and is helpful
to find accuracy bounds when we have to approximate invariant
subspaces of~$\Kaug$, as we explain next.

\section{Non-Invariant Subspaces, Invariance Proximity, and
Approximation Error}

\new{Here, we tackle Problem~\ref{prob:statement}(d)-(e) in Section~\ref{sec:problem-statement}.} 
In the sections above we have provided results connecting the
finite-dimensional invariant subspaces of $\Kaug$ to common invariant
subspaces of the Koopman control family $\{\Kc_{u^*}\}_{u^* \in \Uc}$,
and how these can be used in predicting the evolution of functions on
the common invariant subspace under the trajectories of the control
system. In practice, however, finding exact invariant subspaces that
capture proper information is an arduous task and one might need to
settle for approximations on non-invariant subspaces. In such case, three
fundamental questions immediately arise:
\begin{enumerate}
\item[(Q1)] How can we measure the closeness of a subspace to being
invariant?
\item[(Q2)] How does this measure characterize the approximation error
of the action of the operator on a non-invariant subspace?
\item[(Q3)] How do the previous results regarding the prediction of
observables on the trajectories of the control system extend to the
case of non-invariant subspaces?
\end{enumerate}
These are the questions we tackle in this section. To determine
whether a finite-dimensional subspace $\Sc \subset \Faug$ is invariant
under $\Kaug$ we only need the concept of set inclusion. However, to
quantify how close to invariant a subspace is, we need to be able to
measure angles, lengths, and distances. Therefore, we equip the space
$\Faug$ with an inner product, that induces a norm and, in turn, a
metric\footnote{Even though we aim to approximate a common invariant
  subspace $\Hc \subset \Fc$ under the Koopman control family, our end
  goal is to predict the evolution of observables under the system
  trajectories, i.e., we aim to predict $h(x^+) = h \circ \Tc(x,u)$
  for all $h \in \Hc$ and $(x,u) \in \Xc \times \Uc$. Since
  $h\circ \Tc \in \Faug$, we need to reason with~$\Faug$.}.

\begin{definition}\longthmtitle{Inner Product, Norm, and Metric on
    $\Faug$}\label{def:innerprod-norm-metric}
  An arbitrary inner product\footnote{Since we are working with
    finite-dimensional subspaces, we do not require the inner product
    space $\Faug$ to be complete (Hilbert) or separable.}
  $\innerprod{\cdot}{\cdot}: \Faug \times \Faug \to \cplx$ on $\Faug$
  induces a norm $\| \cdot \|: \Faug \to [0,\infty)$ and a metric
  $\dist: \Faug \times \Faug \to [0,\infty)$ as
  $ \| f \| = \sqrt{\innerprod{f}{f}}$ and $ \dist(f,g) = \|f -
  g\|$. \oprocend
\end{definition}

Since we work with a finite-dimensional subspace that is not
necessarily invariant under the operator, we have to approximate the
action of the operator on the subspace. This approximation is
generally done by performing an orthogonal projection on the subspace,
as explained next.

\begin{definition}\longthmtitle{Linear Predictors on
Finite-Dimensional Subspaces}\label{def:linear-predictors}
Consider the finite-dimensional subspace $\Sc \subset \Faug$ and let
$\Pc_{\Sc}: \Faug \to \Faug$ be the orthogonal projection
operator\footnote{Given an orthonormal basis $\{e_1,\ldots,e_n\}$ for
  $\Sc$, one can calculate the orthogonal projection of $g \in \Faug$
  on $\Sc$ by $\Pc_{\Sc}(g) = \sum_{i=1}^{n} \innerprod{g}{e_i}e_i$.}
on $\Sc$. We define the predictor for the function $\psi \in \Faug$ on
$\Sc$ as \new{$ \Pf_\psi^\Sc := \Pc_{\Sc} \psi$.}  For a vector-valued
function $\Psi = [\psi_1,\ldots,\psi_n]^T$, where $\psi_i \in \Faug$
for $i \in \until{n}$, we define the linear predictor
\new{$\Pf_{ \Psi}^\Sc := [\Pf_{\psi_1}^\Sc, \ldots,
  \Pf_{\psi_n}^\Sc]$}.  We remove the superscript $\Sc$ when it is
clear from the context.
\oprocend
\end{definition}

The properties of the operator $\Pc_{\Sc}$ lead to useful properties
of the linear predictors defined in
Definition~\ref{def:linear-predictors}.

\begin{lemma}\longthmtitle{Properties of Linear
Predictors}\label{l:predictor-properties}
Linear predictors on the finite-dimensional subspace
$\Sc \subset \Faug$ satisfy:
\begin{enumerate}
\item $\Pf_f \in \Sc$ is the best approximation for $f \in \Faug$ on
  $\Sc$, i.e., $\| f - \Pf_f\| \leq \|f - g\|$ for all $g \in \Sc$;
\item $ \Pf_{c_1 f_1 + c_2 f_2} = c_1 \Pf_{f_1} + c_2 \Pf_{f_2}$ for
  all $f_1,f_2 \in \Faug$ and $c_1,c_2 \in \cplx$;
\item let $\Psi$ be a vector-valued function with
  $\Span(\Psi) \subset \Faug$ and let $f = v_f^T \Psi$, where $v_f$ is
  a complex vector of appropriate size. Then,
  $\Pf_f = v_f^T \Pf_\Psi$.  \oprocend
\end{enumerate}
\end{lemma}

The proof of Lemma~\ref{l:predictor-properties} is a direct
consequence of the properties of orthogonal projections and is omitted
for space reasons. 
\new{Lemma~\ref{l:predictor-properties}(a) states that the predictor
  defined in Definition~\ref{def:linear-predictors} is the best
  predictor on the subspace: in this sense, we use the notation
  $f \approx \Pf_f$ when we aim to emphasize that we approximate $f$
  with~$\Pf_f$.  }

We next use the linear predictors to approximate the action of
the operator $\Kaug$ on a non-invariant finite-dimensional subspace
and provide a matrix notation for it.

\begin{lemma}\longthmtitle{Approximating an Operator's Action using
Linear Predictors}\label{l:operator-predictor-matrix}
Any finite-dimensional subspace $\Sc \subset \Faug$ is invariant
under $\Pc_{\Sc}\Kaug$.  Let $\Phi: \Xc \times \Uc \to \cplx^s$ be a
basis for $\Sc$ and let $\tilde{A} \in \cplx^{s \times s}$ be a
matrix such that $\Pc_{\Sc} \Kaug \Phi = \tilde{A} \Phi$. Then,
\begin{enumerate}
\item $\Pf_{\Kaug \Phi} = \tilde{A} \Phi$;
\item for $f \in \Sc$ with description $f = v_f^T \Phi$, where
$v_f \in \cplx^s$, we have $\Pf_{\Kaug f} = v_f^T \tilde{A} \Phi$.
\oprocend
\end{enumerate}
\end{lemma}

Note the parallelism of Lemma~\ref{l:operator-predictor-matrix}
with~\eqref{eq:Koopman-matrix-noninvariant}-\eqref{eq:Koopman-predictor-noninvariant}. Its
proof is a direct consequence of the linearity of $\Kaug$ and
Lemma~\ref{l:predictor-properties}, and is omitted for space reasons.
The prediction error associated with the predictors in
Lemma~\ref{l:operator-predictor-matrix} directly depends on how close
to invariant the space is under the operator $\Kaug$. To capture this,
we define the concept of invariance proximity under an operator.

\begin{definition}\longthmtitle{Invariance
Proximity}\label{def:invariance-proximity}
The \textbf{invariance proximity} of a finite-dimensional subspace
$\Sc \subset \Faug$ under the operator $\Kaug$, denoted
$I_{\Kaug} (\Sc)$, is
\begin{align*}
I_{\Kaug} (\Sc) = \sup_{f \in \Sc, \|\Kaug f\| \neq 0} \frac{ \|
\Kaug f -  \Pf_{\Kaug f}\|}{\| \Kaug f \|}.  \eqoprocend
\end{align*}
\end{definition}

Invariance proximity measures the worst-case relative error of
approximation by projecting the action of $\Kaug$ on~$\Sc$ and
provides an answer to Q2 above.  Invariance proximity does not depend
on the specific basis for the subspace, and is instead a property of
the linear space $\Sc$ and the operator~$\Kaug$.

\begin{proposition}\longthmtitle{Properties of Invariance
Proximity}\label{prop:prop-inv-prox}
Given a finite-dimensional subspace $\Sc \subset \Faug$,
\begin{enumerate}
\item $I_{\Kaug} (\Sc) \in [0,1]$;
\item $I_{\Kaug} (\Sc) = 0$ if~\footnote{The converse also holds if
$\|f-g\|=0$ implies $f=g$ everywhere. This might not hold for
typical norms on function spaces that operate on equivalence
classes and allow for violations of equality on measure-zero
sets.} $\Sc$ is invariant under~$\Kaug$.
\end{enumerate}
\end{proposition}
\begin{pf}
(a) Let $f \in \Faug$ with $\|\Kaug f\| \neq 0$.  Noting that
$\Pf_{\Kaug f} = \Pc_{\Sc} \Kaug f$ is an orthogonal projection on
$\Sc$, we can decompose $\Kaug f$ as $\Kaug f = \Pf_{\Kaug f} + e$,
where $ \langle \Pf_{\Kaug f}, e \rangle = 0$.
Using the definition of the norm induced by the inner product then
yields $ \|\Kaug f\|^2 = \|\Pf_{\Kaug f}\|^2 + \|e\|^2$.  Therefore,
$\|e\| \leq \|\Kaug f\|$ and we can write 
\begin{align*}
\frac{ \| \Kaug f -  \Pf_{\Kaug f}\|}{\| \Kaug f \|} = \frac{ \| e
\|}{\| \Kaug f \|} \leq 1. 
\end{align*}
Since this inequality holds for all functions $f \in \Faug$ where
$\| \Kaug f \| \neq 0$, we deduce $I_{\Kaug} (\Sc) \leq
1$. Moreover, by definition of $I_{\Kaug} (\Sc)$ and the fact that
norms are nonnegative, we conclude $I_{\Kaug} (\Sc) \geq 0$,
completing the proof.

(b) If $\Sc$ is invariant under $\Kaug$, we have $\Kaug f \in \Sc$
and therefore $\| \Kaug f - \Pf_{\Kaug f}\| =0$ for all $f \in
\Sc$. Hence, $I_{\Kaug} (\Sc) = 0$. \qed
\end{pf}

Proposition~\ref{prop:prop-inv-prox} means that invariance
proximity provides an answer to Q1 above.  The next result extends to
non-invariant subspaces the results on prediction of the evolution of
functions in $\Fc$ under the control system~\eqref{eq:control-system},
providing an answer to~Q3.

\begin{theorem}\longthmtitle{Approximate Input-State Separable Form
and Accuracy
Bound}\label{t:invariance-proximity-bounds-error-in-F}
Let $\Sc \subset \Faug$ be a finite-dimensional normal subspace
comprised of input-state separable combinations. Let
$\Phi(x,u)= G(u) H(x)$ be a decomposition of a basis for $\Sc$ where
$G: \Uc \to \cplx^{s \times l}$ and $H: \Xc \to \cplx^{l}$ for some
$l \leq s$ (here, $s \in \naturals$ is the dimension of $\Sc$). Let
$H_e$ and $h_e$ be the control-independent extensions of $H$ and
$h \in \Span(H)$, resp. Then,
\begin{enumerate}
\item $\Pf_{h \circ \Tc} = \Pf_{\Kaug h_e}$ for all
$h \in \Span(H)$, and $\Pf_{H \circ \Tc} = \Pf_{\Kaug H_e}$;
\item without loss of generally, assume $\Phi$ is in normal form,
i.e., $G(u) = I_{\Uc}^{l \times l}(u)$ if $l = s$ or
$G(u) = [I_{\Uc}^{l \times l}(u)^T, \tilde{G}(u)^T]^T$ if $s
>l$. Moreover, let $\tilde{A} \in \cplx^{s \times s}$ be a matrix
such that $\Pc_{\Sc} \Kaug \Phi = \tilde{A} \Phi$ (note that
$\tilde{A}$ exists because $\Sc$ is invariant under
$\Pc_{\Sc}\Kaug$). If $s > l$, consider the block-decomposition of
$\tilde{A}$,
\begin{align*}
\tilde{A} =
\begin{bmatrix}
\tilde{A}_{11} & \tilde{A}_{12}
\\
\tilde{A}_{21} & \tilde{A}_{22}
\end{bmatrix},  
\end{align*}
where $\tilde{A}_{11} \in \cplx^{l \times l}$,
$\tilde{A}_{12} \in \cplx^{l \times (s-l)}$,
$\tilde{A}_{21} \in \cplx^{(s-l) \times l}$, and
$\tilde{A}_{22} \in \cplx^{(s-1) \times (s-l)}$. Then, the associated
approximate input-state separable dynamics can be written as
\begin{align}\label{eq:approx-separable-form}
H(x^+) = H \circ \Tc(x,u) \approx \Pf_{H \circ \Tc} (x,u)= \Ac(u) H(x),
\end{align}
where,  for each $u \in \Uc$,
\begin{align*}
&\Ac(u) = \tilde{A}_{11} + \tilde{A}_{12} \tilde{G}(u),
&&\text{if $s > l$}. 
\\ 
&\Ac(u) = \tilde{A} \, I_{\Uc}^{l \times l}(u) = \tilde{A},
&&\text{if $s = l$}.  
\end{align*}

\item for all $h \in \Span(H)$ with description $h = v_h^T H$,
$v_h \in \cplx^l$,
\begin{align*}
h(x^+) = h \circ \Tc(x,u) \approx  \Pf_{h \circ \Tc} (x,u) =
v_h^T \Ac(u) H(x); 
\end{align*}

\item for all $h \in \Span(H)$ with $ \|h \circ \Tc\| \neq 0$, the
  predictor's relative error is bounded by the invariance proximity of
  $\Sc$ under $\Kaug$, \new{as described by the inequality}
\begin{align*}
\frac{\left\| h \circ \Tc - \Pf_{h \circ \Tc} \right\|}{\| h
\circ \Tc \|} \leq I_{\Kaug}(\Sc) .
\end{align*}
\end{enumerate}
\end{theorem}
\begin{pf}
(a) By Definition~\ref{def:linear-predictors},
$\Pf_{h \circ \Tc} = \Pc_{\Sc} (h \circ \Tc)$. Using
Lemma~\ref{l:Kaug-captures-function-evolutions-in-F}, we have
$h \circ \Tc = \Kaug h_e$. Hence,
$\Pf_{h \circ \Tc} = \Pc_{\Sc} \Kaug h_e = \Pf_{\Kaug h_e}$. The
statement regarding $H$ follows directly by applying this to each
element of the equality $\Pf_{H \circ \Tc} = \Pf_{\Kaug H_e}$.

(b) We need to prove the rightmost equality
in~\eqref{eq:approx-separable-form}, since the rest follow directly
from their definitions. From part~(a), and using the vector-valued
notation in Remark~\ref{r:notation-overload}, we have
\begin{align}\label{eq:predictor-He}
\Pf_{H \circ  \Tc} = \Pf_{\Kaug H_e} = \Pc_{\Sc} \Kaug H_e.
\end{align}
For the case $s = l$, we use
Lemma~\ref{l:control-indep-properties}(b) to write
$\Phi(x,u) = I_{\Uc}^{l \times l}(u) H(x) = H_e(x,u)$. Hence, noting
that $\Pc_{\Sc} \Kaug \Phi = \tilde{A} \Phi$, we have
$\Pc_{\Sc} \Kaug H_e (x,u) = \tilde{A} H_e(x,u) = \tilde{A}
I_{\Uc}^{l \times l}(u) H(x)$. Using~\eqref{eq:predictor-He}, we can
write
$\Pf_{H \circ \Tc}(x,u) = \tilde{A} H_e (x,u)= \tilde{A} I_{\Uc}^{l
\times l}(u) H(x)$, which completes the proof.

Next, we turn our attention to the case $s>l$. Using
Lemma~\ref{l:control-indep-properties}(b), one can write
\begin{align}\label{eq:phi-decomposition}
\Phi(x,u) = [H_e(x,u)^T, (\tilde{G}(u) H(x))^T]^T.
\end{align}
Multiplying both sides of $\Pc_{\Sc} \Kaug \Phi = \tilde{A} \Phi$
from the left by $W = [I_{l \times l}, 0_{l \times (s-l)}]$, and
using~\eqref{eq:phi-decomposition}, the decomposition of
$\tilde{A}$, the properties of the vector-valued notation in
Remark~\ref{r:notation-overload}, and the linearity of the operator
$\Pc_{\Sc} \Kaug$, one can write
\begin{align*}
\Pc_{\Sc} \Kaug H_e
&= W \Pc_{\Sc} \Kaug \Phi = W \tilde{A} \Phi 
\\
&= (\tilde{A}_{11} I_{\Uc}^{l \times l}(u) + \tilde{A}_{12}\tilde{G}(u)) H(x).
\end{align*}
The statement then follows from equation~\eqref{eq:predictor-He} and
the fact that $I_{\Uc}^{l \times l}(u) = I$ for all $u \in \Uc$.

(c) We need to prove the rightmost equality
$\Pf_{h \circ \Tc}(x,u) = v_h^T\Ac(u) H(x)$, since the rest follow
directly from their definitions.  By hypothesis
$h \circ \Tc = v_h^T H \circ \Tc$; hence, from
Lemma~\ref{l:predictor-properties}(c), we have
$\Pf_{h \circ \Tc} = v_h^T \Pf_{H \circ \Tc}$. The result then
follows from~\eqref{eq:approx-separable-form}.

(d) By Proposition~\ref{p:normal-space-u-independent}, and using the
definition of invariance proximity, for all $h \in \Span(H)$ with
$\|h \circ \Tc\| \neq 0$,
\begin{align*}
\frac{\| \Kaug h_e - \Pf_{\Kaug h_e}\|}{\| \Kaug h_e \|} \leq I_{\Kaug}(\Sc).
\end{align*}
The statement then follows from the fact that
$\Kaug h_e = h \circ \Tc $
(cf.~Lemma~\ref{l:Kaug-captures-function-evolutions-in-F}) and
part~(a). \qed
\end{pf}

This result can be viewed as an analog of
Theorem~\ref{t:normal-space-common-invariant} for non-invariant
subspaces. Theorem~\ref{t:invariance-proximity-bounds-error-in-F}
allows to approximate models in the input-state separable form
(cf. Theorem~\ref{t:common-invariant-separates}) by approximating a
single normal invariant subspace of $\Kaug$, which is significantly
easier than working with the KCF directly. Moreover, the concept of
invariance proximity provides a bound for approximation errors on the
entire subspace.  This has important implications for the validity and
approximation accuracy of common Koopman-inspired descriptions of the
control system~\eqref{eq:control-system},
cf. Lemmas~\ref{l:switched-linear-separable}
and~\ref{l:linear-bilinear-separable}.

\section{Implications for Robust Data-driven Learning}
In this section we illustrate how the results of the paper can be used
in data-driven modeling of control systems \new{and answer Problem~\ref{prob:statement}(f) in Section~\ref{sec:problem-statement}}. We provide an algorithmic
description that specifies how to process the data, the choice of
inner product space, and the formulation for the dictionary learning.

\subsection{Gathering Data for the Augmented Koopman Operator}
Our strategy for learning relies on using
Theorem~\ref{t:invariance-proximity-bounds-error-in-F} to
approximate an input-state separable form and bound the prediction
error for all functions in the identified subspace. This result
employs the augmented Koopman operator associated with the augmented
system~\eqref{eq:augmented-system} and, instead, we can only collect
trajectory data from the original control
system~\eqref{eq:control-system}. This mismatch can be easily
reconciled as we explain next.

Let $\{x_i\}_{i=1}^N \subset \Xc$, $\{u_i\}_{i=1}^N \subset \Uc$, and
$\{x_i^+\}_{i=1}^N \subset \Xc$ be state and input data from
trajectories of system~\eqref{eq:control-system} such that
\begin{align}\label{eq:data-snapshots-control}
x_i^+ = \Tc(x_i, u_i), \quad \forall i \in \until{N}.
\end{align}
A close look at the definition of $\Taug$
in~\eqref{eq:augmented-system} reveals that it does not alter the
input signal, i.e., if we apply it on the state-input pair $x_i, u_i$
for all $i \in \until{N}$, we get
$\Taug(x_i,u_i) = (\Tc(x_i,u_i),u_i) = (x_i^+, u_i)$.  Therefore, we
already have access to all the information $\Taug$ generates: the
first element returned by $\Taug$ is exactly the action of the control
system $\Tc$ that we have measured
in~\eqref{eq:data-snapshots-control} and the second element is exactly
the input (without any change) to $\Tc$, again measured
in~\eqref{eq:data-snapshots-control}. For convenience, we gather these
data snapshots for $\Taug$ in snapshot matrices as follows
\begin{align}\label{eq:data-augmented-system}
&X = [x_1, \ldots, x_N] \in \real^{n \times N},
&&
X^+ =[x_1^+,
\ldots,
x_N^+]\in
\real^{n \times
N}, 
\nonumber
\\
&U = [u_1, \ldots, u_N]\in \real^{m \times N},  &&U^+ = U.
\end{align}
Note that even though matrix $U^+$ does not capture additional
information we have created it, since it is a part of the
corresponding state for $\Taug$. To apply existing numerical methods
such as EDMD on $\Kaug$, we gather the augmented state snapshots of
$\Taug$ as
\begin{align}\label{eq:state-snapshots-augmented-system}
&Z = [X^T, U^T]^T \in \real^{(n+m) \times N},  
\nonumber \\
&Z^+ = [(X^+)^T, (U^+)^T]^T \in \real^{(n+m) \times N}.
\end{align}

\subsection{Choice of Inner Product Space}
The results in the previous sections can be used for subspace learning
on any arbitrary inner product space. Here we focus on the most
popular inner product space in the literature that is used for the
EDMD method~\cite{MOW-IGK-CWR:15,MK-IM:18}. Consider the empirical
measure $\mu_Z$ defined by
\begin{align}\label{eq:empirical-measure-z}
\mu_Z = \frac{1}{N} \sum_{i=1}^{N} \delta_{z_i},
\end{align}
where $\delta_{z_i}$ is the Dirac measure at point $z_i$, the $i$th
column of matrix $Z$ defined
in~\eqref{eq:state-snapshots-augmented-system}. We then choose the
space $L_2(\mu_Z)$ comprised of functions on the domain
$\Xc \times \Uc$.  Under this choice, given any basis
$\Phi: \Xc \times \Uc \to \real^s$ with \emph{real-valued} elements
(cf. Remark~\ref{r:real-dictionary}) for the finite-dimensional (with
dimension $s$) normal subspace $\Sc$, the matrix $\tilde{A}$ in the
hypotheses of Theorem~\ref{t:invariance-proximity-bounds-error-in-F}
is the EDMD solution applied on dictionary $\Phi$ and data
in~\eqref{eq:state-snapshots-augmented-system} (cf.
Section~\ref{subsec:EDMD}), i.e.,
\begin{align}\label{eq:data-driven-EDMD}
\tilde{A} = \Phi(Z^+) \Phi(Z)^ \dagger.
\end{align}
Moreover, under the condition that $\Phi(Z)$ and $\Phi(Z^+)$ have full
row rank, the invariance proximity turns into the square root of the
consistency index (cf.  Section~\ref{subsec:consistency-index}) and
has the following closed-form expression
\begin{align}\label{eq:data-driven-invariance-proximity}
I_{\Kaug}(\Sc)
&= \sqrt{\ic(\Phi,Z,Z^+)} 
\nonumber \\
&= \sqrt{\lambda_{\max}\big(I - \Phi(Z^+) \Phi(Z)^ \dagger \Phi(Z)
\Phi(Z^+)^\dagger\big)} .
\end{align}
We use~\eqref{eq:data-driven-invariance-proximity} to formulate an
optimization-based learning problem for modeling the control system.

\subsection{Optimization-Based Subspace Learning}
Based on Theorem~\ref{t:invariance-proximity-bounds-error-in-F}(d),
the invariance proximity determines the accuracy of the model
provided on a given normal subspace. Hence, we formulate an
optimization problem to find an accurate model by minimizing the
invariance proximity over a parametric family of normal spaces with
basis $\Phi$ in normal form~\eqref{eq:normal-form-dictionary} as
\begin{align}\label{eq:optimization}
&\underset{ \Phi \in \text{PF}}{\text{minimize}} \; I_{\Kaug}(\Sc)
\! & \Leftrightarrow \!
&&\underset{ \Phi \in
\text{PF}}{\text{minimize}} \;
\sqrt{\ic(\Phi,Z,Z^+)} \;, 
\end{align}
where $\text{PF}$ is the parametric family of choice (e.g.,~neural
networks, polynomials), $\Sc = \Span(\Phi)$, and one can use the
closed-form solution of the invariance proximity
in~\eqref{eq:data-driven-invariance-proximity}.  Note that depending
on the choice of the parametric family, the optimization
problem~\eqref{eq:optimization} is generally non-convex.

We make the following observations regarding the optimization
problem~\eqref{eq:optimization} and its properties:

\paragraph*{Alternative formulation for efficiency and numerical
resiliency to finite-precision errors}
Using the closed-form expression for invariance proximity
in~\eqref{eq:data-driven-invariance-proximity} requires calculating
the maximum eigenvalue of
$M_C = I - \Phi(Z^+) \Phi(Z)^ \dagger \Phi(Z)
\Phi(Z^+)^\dagger$. This matrix has spectrum in $[0,1]$,
cf.~\cite[Lemma~1]{MH-JC:23-csl}. Many software packages for finding
maximum eigenvalues rely on iterative methods that are sensitive to
the separation between the largest and second largest
eigenvalues. To avoid numerical issues, one can use $\Tr(M_C)$
instead of $\lambda_{\max}(M_C)$, as justified by
\begin{align*}
\frac{1}{s} \Tr(M_C) \leq \lambda_{\max}(M_C) \leq \Tr(M_C),
\end{align*}
where $s$ is the dimension of $M_C$. Note that the inequalities follow
from the fact that the spectrum of $M_C$ belongs to~$[0,1]$.

\paragraph*{Equivalence to robust minimax problem}
Based on Theorem~\ref{t:RRMSE-bound-sprad-consistency}, the
optimization problem~\eqref{eq:optimization} is equivalent to the
following robust minimax problem
\begin{align*}
\underset{ \Phi \in \text{PF}}{\text{minimize}} \max_{f \in \Sc,
\| \Kaug f  \|_{L_2(\mu_Z)} \neq 0} \frac{\| \Kaug f - 
\Pf_{\Kaug f} \|_{L_2(\mu_Z)}}{\| \Kaug f  \|_{L_2(\mu_Z)}},
\end{align*}
where $\Sc = \Span(\Phi)$ and $\mu_Z$ is defined
in~\eqref{eq:empirical-measure-z}. This equivalence makes it clear
that optimization~\eqref{eq:optimization} minimizes the worst-case
error on the subspace, does not depend on the choice of basis, and
is not sensitive to the scaling of variables.

For the readers' convenience, Algorithm~\ref{algo:data-driven-model}
summarizes the steps described above to learn input-state separable
models.

\new{
\begin{algorithm}
\caption{Learning Input-State Separable Models} \label{algo:data-driven-model}
\begin{algorithmic}[1] 
\Statex $\triangledown$ \textbf{Data acquisition}
\State Gather data according
to~\eqref{eq:data-augmented-system}
\State Form matrices $Z$ and $Z^+$ according to~\eqref{eq:state-snapshots-augmented-system}
\vspace*{5pt}
\Statex $\triangledown$ \textbf{Approximate the action of  $\Kaug$ on normal space}
\State Choose parametric family  of normal dictionaries
(e.g., neural networks, polynomials)
with real-valued elements in the form~\eqref{eq:normal-form-dictionary}
\State Obtain $\Phi^*$  by solving~\eqref{eq:optimization} 
\State Calculate
$\tilde{A} = \Phi^*(Z^+) \Phi^*(Z)^ \dagger$
\vspace*{5pt}
\Statex $\triangledown$ \textbf{Determine input-state separable model}
\State Find  input-state separable form via
Theorem~\ref{t:invariance-proximity-bounds-error-in-F}(b)
\end{algorithmic}
\end{algorithm}
}

\begin{example}\longthmtitle{DC Motor with Nonlinear Multiplicative
    Input Injection} {\rm Consider the DC motor\footnote{This system
      is a modified version of the experimental
      study~\cite{SDB-HU:98}, used as an example
      in~\cite{MK-IM-automatica:18}. We have added the nonlinear
      function $f:\Uc \to \real$ in the mechanism generating the field
      current.}
\begin{align}\label{eq:DC-motor}
\dot{x}_1
&= -(R_a /L_a) x_1 - (k_m/ L_a) x_2 f(u) + u_a/L_a, 
\nonumber
\\
\dot{x}_2
&= -(B /J) x_2 + (k_m/ J) x_1 f(u) -  \tau_l/J,
\end{align}
where $x_1$ is the armature current, $x_2$ is the angular velocity,
$x = [x_1, x_2]^T$ is the state vector and $u$ is the input.  The
value of the parameters are $R_a = 12.345$, $L_a = 0.314$,
$k_m = 0.253$, $u_a =60$, $B=0.00732$, $\tau_l = 1.47$, and
$J = 0.00441$. We consider two choices for $f$:
(i)~$f(u) = 2 \tanh(u)$ (saturated input) and
(ii)~$f(u) = 2 \tanh(u \cos(u))$ (saturated non-monotone input).
\new{The normal operating range of the motor and the input set are
  $\Mc = [-5,15] \times [-250,125]$ and $\Uc = [-4,4]$
  resp. We define the system's state space $\Xc$ as the
  reachable set from $\Mc$ given all possible input signals taken from
  $\Uc$.

\textit{Data:} We run $10^4$ experiments with constant inputs
and length $50 \, ms$ with uniformly selected initial conditions
from  $\Mc$ and
inputs from $\Uc = [-4,4]$. We sample the trajectories with time step
$\Delta t = 5 \, ms$, resulting in a total of $10^5$ data
snapshots. Out of this data set, we select half as the training
data set and the rest as the test data set.
}

\textit{Parametric Families for Comparison:} Our aim is to
compare the effectiveness of our methods with widely used
lifted-linear (an extension of DMD with
control~\cite{JLP-SLB-JNK:16}) and bilinear Koopman-inspired
forms as
\begin{align*}
\Psi_l(x^+) &= A_l \Psi_l(x) + B_l u,
\\
\Psi_b(x^+) &= A_b \Psi_b(x) + B_b \Psi_b(x)u,
\end{align*}
We use data to learn all the models with dimension four. For the
input-state separable model\footnote{\new{The choice of dimensions is
    problem dependent. In the normal
    basis~\eqref{eq:normal-form-dictionary}, $l$ determines the
    dimension of input-state separable model. The choice of parameter
    $s$ determines how rich the input-dependent part is. We recommend
    the user to set $s \geq 2 \, l$ so that the  matrix-valued function
    $\tilde{G}(u)$ in~\eqref{eq:normal-form-dictionary} has more (or equal) rows than columns. Moreover,
    for $l$, we recommend starting from a small number and increasing
    if the invariance proximity of the trained model is larger than
    the desired level.}}, in the normal
basis~\eqref{eq:normal-form-dictionary}, we set the dimension of
normal space $s = 20$ and the dimension of the input-state separable
model as $l = 4$. We model the functions $H(x)$ and $\tilde{G}(u)$
in~\eqref{eq:normal-form-dictionary} by two residual neural
networks~\citep{KH-XZ-SR-JS:16} comprised of 5 residual blocks each
with 64 neuron per hidden layer and ReLU activation functions. We also
fix the first two elements of $H(x)$ to be the state vector
corresponding to the system. For the linear and bilinear models, we
set the functions $\Psi_l(x)$ and $\Psi_b(x)$ to be the same type of
neural network used for $H(x)$ in the input-state separable model. To
learn the input-state separable model we use
Algorithm~\ref{algo:data-driven-model}. We train the neural networks
for lifted linear and bilinear models by minimizing the following
typical least norm residual errors
\begin{align}\label{eq:linear-bilinear-learning}
&\underset{ \Psi_l \in \text{PF}}{\text{minimize}} \quad
\|\Psi_l(X^+) - A_l^* \Psi_l(X) - B_l^* U \|_F,  
\nonumber
\\ 
&\underset{ \Psi_b \in \text{PF}}{\text{minimize}} \quad
\|\Psi_b(X^+) - A_b^* \Psi_b(X) - B_b^* \Psi_b(X) \! \cdot
\! U \|_F, 
\end{align}
where $\Psi_b(X)\!  \cdot \! U$ denotes the column-wise product
of $\Psi_b(X)$ and $U$. Moreover, $A_l^*$ and $B_l^*$ are the
best parameters by minimizing the same cost function over $A_l$
and $B_l$ instead of $\Psi_l$
(see e.g.,~\citep[Section~4]{MK-IM-automatica:18}). $A_b^*$ and
$B_b^*$ are computed similarly. Note that the optimization
problems in~\eqref{eq:linear-bilinear-learning} are solved in
$\Psi_l$ and $\Psi_b$ and are nonconvex.

\textit{Training and Practical Considerations:} We randomly initialize
the neural networks. The networks for $H(x)$ in the input state
separable model, and $\Psi_l$ and $\Psi_b$ start from the same initial
weights and biases. To make sure all variables are in the same scale
and one does not dominate the others, we do a change of coordinates by
$x_1 \mapsto 0.1 \, x_1$, $x_2 \mapsto 0.004 \, x_2$ and
$u \mapsto 0.25 \, u$. We scale back to the original coordinates after
training. One also might benefit from regularization in case of
overfitting. However, note that the robust minimax problem used in our
method is resilient to overfitting since it considers all uncountably
many functions in the vector space and we have not used any
regularization. To train the neural networks, we use the Adam
method~\cite{DPK-JB:15} with batch size of 200. We train the networks
for 500 epochs while decreasing the learning rate linearly from
$5 \times 10^{-4}$ to $10^{-6}$.  Finally, we use the formula in
Theorem~\ref{t:invariance-proximity-bounds-error-in-F}(b) to build an
input-state separable model based on the augmented operator.

\textit{Evaluation and Comparison:} To evaluate the accuracy of
models, we create a piecewise constant random input with time
step $\Delta t$ for $600$ time steps (or 3 seconds) and compare
the learned models' response with the actual system trajectories
generated by~\eqref{eq:DC-motor}. Figures~\ref{fig:random-input}
shows the generated input signal used for
comparison. Figure~\ref{fig:DC-motor-saturated} shows the
angular velocity of the motor with nonlinear input injection
$f(u) =2 \tanh(u)$ compared to predictions derived by the
input-state separable, lifted linear, and lifted bilinear
models. Figure~\ref{fig:DC-motor-periodic} depicts the same
comparison for the case where $f(u) = 2 \tanh(u \cos(u))$ in
system~\eqref{eq:DC-motor}. Clearly, in both cases
$f(u) = 2 \tanh(u)$ and $f(u) = 2 \tanh(u \cos(u))$ the
input-state separable model outperforms the other
methods. Moreover, by comparing
Figure~\ref{fig:DC-motor-saturated} with
Figure~\ref{fig:DC-motor-periodic}, one can see that the more
nonlinear the system is in the input, the less accurate the
lifted linear and bilinear models become (even for short-term
predictions).
\oprocend

\begin{figure}[htb]
\centering 
{\includegraphics[width=.50\linewidth]{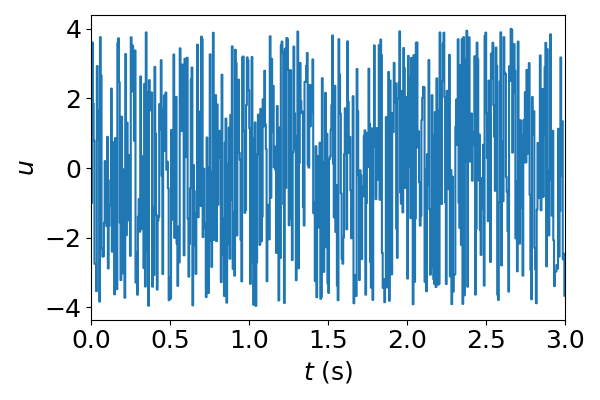}}
\caption{The piecewise constant random signal used for comparing different models with $\Delta t = 5 \, ms$ and length of 600 steps.}\label{fig:random-input}
\vspace*{-1ex}
\end{figure}

\begin{figure}[htb]
\centering 
{\includegraphics[width=.49\linewidth]{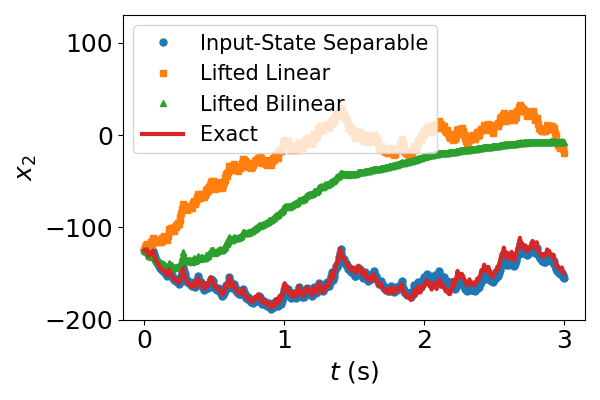}}
{\includegraphics[width=.49\linewidth]{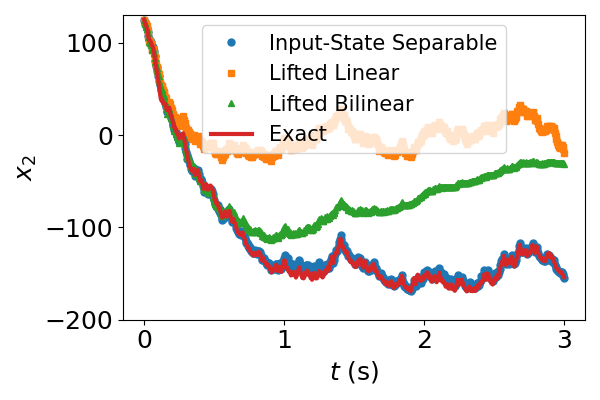}}
\caption{The angular velocity, $x_2$, of the DC motor
in~\eqref{eq:DC-motor} with $f(u) = 2 \tanh(u)$ and
predictions by input-state separable (our method), lifted
linear, and lifted bilinear models. The trajectories start
from two initial conditions $x_0 = [0, -125]$ (left), and
$x_0 = [0, 125]$ (right).}\label{fig:DC-motor-saturated}
\vspace*{-1ex}
\end{figure}

\begin{figure}[htb]
\centering 
{\includegraphics[width=.48\linewidth]{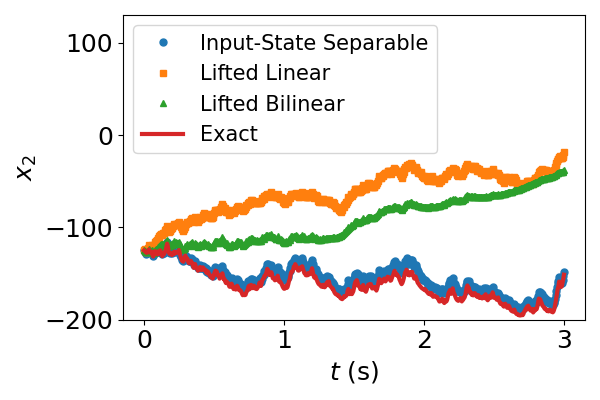}}
{\includegraphics[width=.49\linewidth]{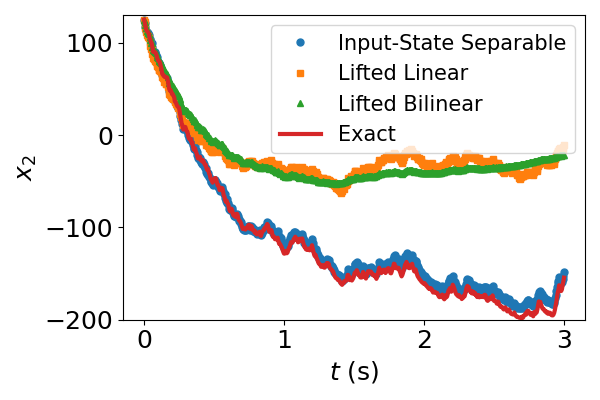}}
\caption{The angular velocity, $x_2$, of the DC motor in~\eqref{eq:DC-motor} with $f(u) = 2 \tanh(u \cos(u))$ and predictions by input-state separable (our method), lifted linear, and lifted bilinear models. The trajectories start from two initial conditions  $x_0 = [0, -125]$ (left), and $x_0 = [0, 125]$ (right).}\label{fig:DC-motor-periodic}
\vspace*{-1ex}
\end{figure}
%

 \textit{Implications for Control Design:} Here, we demonstrate the
 implications of using our framework as enabler for control
 applications via a simple Model Predictive Control (MPC)
 example. Consider the following receding horizon control problem for a
 given time-step $k \in \naturals$
 \begin{align}\label{eq:vanilla-mpc}
   &\underset{\hat{\mathbf{x}}_{k+1}, \mathbf{u}_k}{\text{minimize}}
   && J (\hat{\mathbf{x}}_{k+1}, \mathbf{u}_k)
      \nonumber
   \\
   & \text{subject to:}
   && u_{k+i} \in \Uc, \; \forall i \in \{0,\ldots, H-1\},
      \nonumber
   \\
   & && z_{k+i+1} = \operatorname{model}(z_{k+i}, u_{k+i}), \; \forall
        i \in \{0,\ldots, H-1\}, 
        \nonumber
   \\ 
   & && \hat{x}_{k+i+1} = C z_{k+i+1},  \; \forall i \in \{0,\ldots, H-1\},
        \nonumber
   \\
   & && z_k = \Psi_{\operatorname{model}}(x_k),
 \end{align}
 where
 $\hat{\mathbf{x}}_{k+1} := (\hat{x}_{k+1},\ldots, \hat{x}_{k+H})$,
 $\mathbf{u}_k :=(u_k, \ldots, u_{k+H-1})$,
 $J (\hat{\mathbf{x}}_k, \mathbf{u}_k) : = \sum_{i = 0}^{H-1}
 [(\hat{x}_{k+i+1} - r_{k+i+1})^T Q (\hat{x}_{k+i+1} - r_{k+i+1}) + R
 \|u_{k+i}\|^2]$ with $(r_{k+i})_{i=1}^H$ being the reference
 trajectory we aim to track and $H$ is the horizon length. Moreover,
 in~\eqref{eq:vanilla-mpc}, ``$\operatorname{model}$'' represents one of
 the three models in this section; namely, input-state separable,
 lifted bilinear, or lifted linear. In addition,
 $\Psi_{\operatorname{model}}$ is the appropriate ``lifting function''
 of the state corresponding to the chosen model (for input state
 separable model $\Psi_{\operatorname{model}} (x)= H(x)$) and 
 $C$ is a matrix that extract the elements of the state from
 $\Psi_{\operatorname{model}}$. In this example, we use the following
 parameters
 \begin{align*}
   &C =
     \begin{bmatrix}
       1 & 0 & 0 & 0 \\ 0 & 1 & 0 & 0
     \end{bmatrix},
   &&Q =
      \begin{bmatrix}
        0 & 0
        \\
        0 & 1
      \end{bmatrix},
   &&&R = 0.1,
   &&&&
        H = 20, 
 \end{align*}
 and run the MPC scheme iteratively\footnote{Since all the models used
   here have closed-form solutions given an initial condition and an
   input sequence, we used the closed-form solution to directly
   incorporate the dynamic constraints in~\eqref{eq:vanilla-mpc} into
   the cost function resulting in an optimization problem in variable
   $\mathbf{u}_k$ which we have solved with limited-memory
   Broyden–Fletcher–Goldfarb–Shanno (L-BFGS) algorithm with box
   constraints~\cite{RHB-PL-JN-CZ:95,CZ-RHB-PL-JN:97}.} at each time
 step $k \in \naturals_0$ and apply the first step of the computed
 input to make sure the motor's angular velocity tracks a periodic
 reference signal. At the first time step ($k=0$), the initial
 condition for the optimization~\eqref{eq:vanilla-mpc} is randomly
 selected. After the first time step ($k \geq 1$), we use the
 previously calculated optimal input sequence to generate an initial
 condition for the current step. The procedure is as follows: the last
 element of initial condition $\mathbf{u}_k$ for
 optimization~\eqref{eq:vanilla-mpc} is randomly selected, while the
 first $H-1$ elements are chosen to be the last $H-1$ elements of the
 optimal input sequence from the previous step ($\mathbf{u}_{k-1}$).

 Figure~\ref{fig:MPC-comparison} shows the comparison between MPC
 schemes via different Koopman based models for
 system~\eqref{eq:DC-motor} with $f(u) =2 \tanh(u)$ and
 $f(u) = 2 \tanh(u \cos(u))$. As Figure~\ref{fig:MPC-comparison}
 clearly depicts, there is a significant difference in the performance
 of controller for different types of nonlinearity in the input. For
 the simple case where $f(u) = 2 \tanh(u)$, the lifted bilinear and
 input-state separable models give close performance since, if $u$ is
 reasonably away from the edges of the interval $\Uc = [-4,4]$, the
 system~\eqref{eq:DC-motor} can be accurately approximated by a control
 affine system and thus bilinear models are appropriate for it (see
 e.g.~\cite{DG-DAP:21}). However, for the non-monotonic nonlinearity,
 $f(u) = 2 \tanh(u \cos(u))$, in~\eqref{eq:DC-motor}, the approximation
 by a control-affine system is no longer valid and the lifted bilinear form
 is no longer able to capture the system's behavior even with an MPC
 controller in the loop.

 \begin{figure}[htb]
   \centering 
   {\includegraphics[width=.49\linewidth]{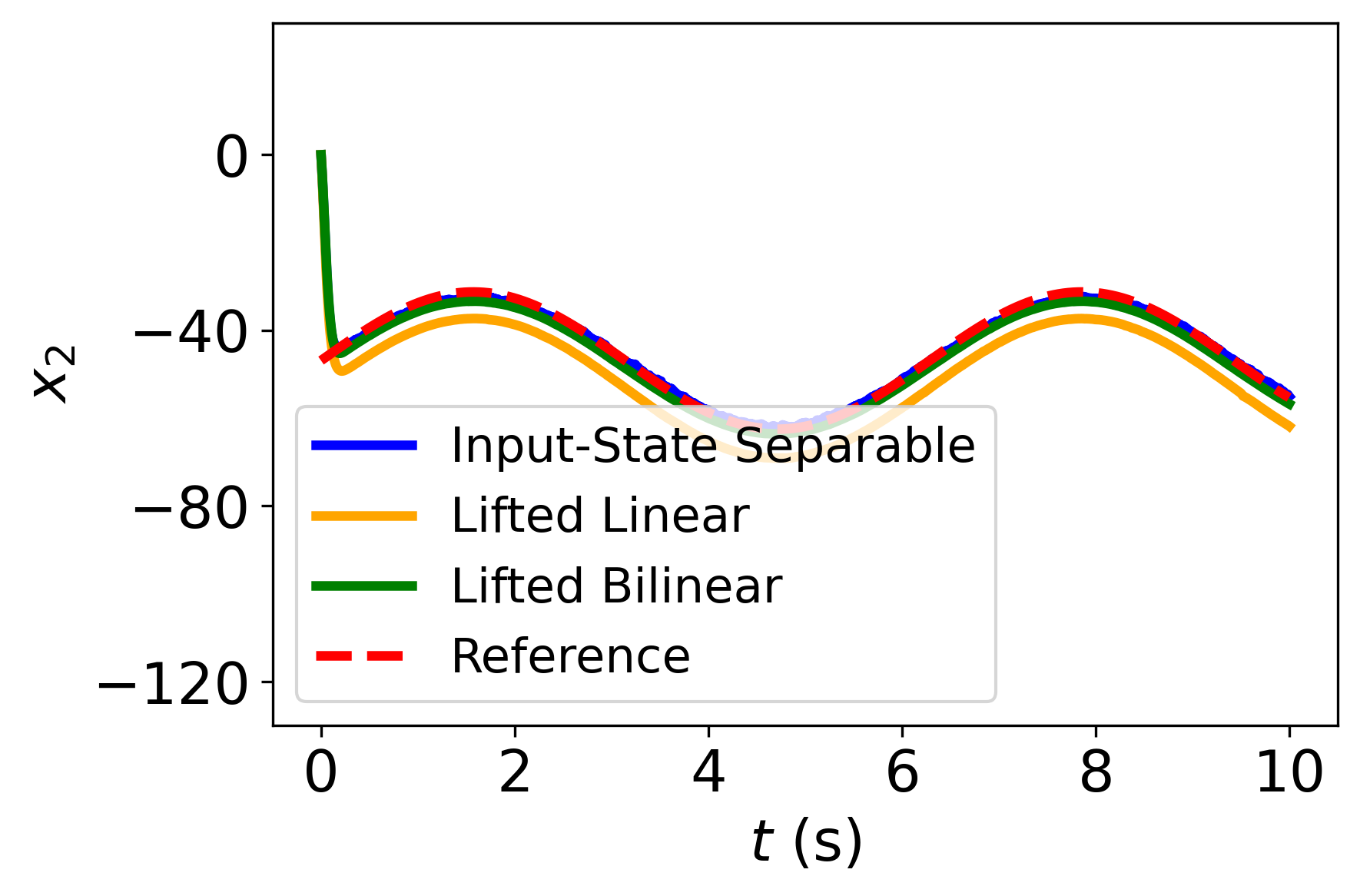}}
   {\includegraphics[width=.49\linewidth]{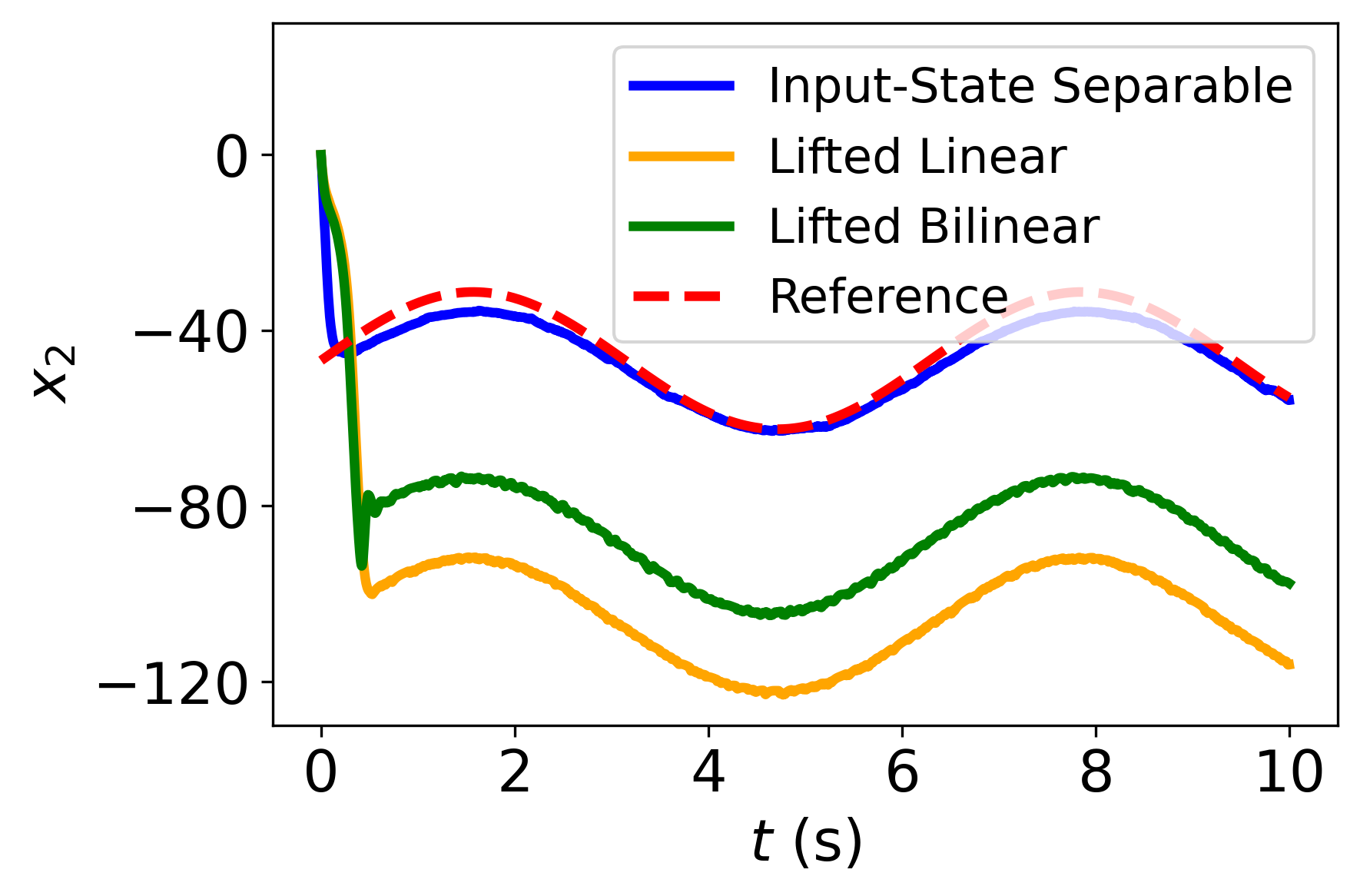}}
   \caption{The angular velocity, $x_2$, of the DC
     motor~\eqref{eq:DC-motor} under MPC control vs periodic reference
     signal for input-state separable, lifted linear, and lifted
     bilinear models. The plot on the left corresponds to the DC motor
     with saturated input nonlinearity ($f(u) = 2 \tanh(u)$
     in~\eqref{eq:DC-motor}) and the plot on the right corresponds to
     the DC motor with saturated non-monotone input nonlinearity
     ($f(u) = 2 \tanh(u \cos(u))$
     in~\eqref{eq:DC-motor}).}\label{fig:MPC-comparison}
   \vspace*{-1ex}
 \end{figure}
}
\end{example}

\section{Conclusions}
We have presented the notion of Koopman Control Family (KCF), a
theoretical framework for modeling general nonlinear control
systems. We have shown that the KCF can fully characterize the
behavior of a control system on a (potentially infinite-dimensional)
function space. To build finite-dimensional models, we have introduced
a generalized notion of subspace invariance, leading to a universal
finite-dimensional form which we refer to as \emph{input-state
separable}. Remarkably, the commonly-used lifted linear, bilinear,
and switched linear models are all special cases of the input-state
separable form. We have provided a complete theoretical analysis
accompanied by discussions on usage in data-driven
applications. Future work will build on the results of the paper
to develop strategies for control design, such as using the
closed-form solution of the input-state separable models to provide
computational gains and performance guarantees for model predictive
control as well as extending switching-based linear control designs
to the case of uncountable and potentially unbounded input sets.  We
also aim to build on the proposed framework to determine reachable
and control-invariant sets. We also aim to explore additional
structures that the KCF might enjoy for special classes of nonlinear
systems such as control-affine and monotone systems.

{
\small
\bibliographystyle{plainnat}
\bibliography{alias,JC,Main,Main-add}
}

\medskip

\setlength{\intextsep}{0pt}%
\setlength{\columnsep}{5pt}%
\begin{wrapfigure}{l}{0.25\linewidth}
\centering
\includegraphics[width=\linewidth]{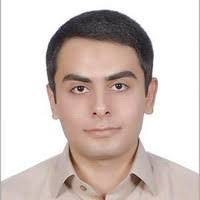}
\end{wrapfigure}
\textbf{Masih Haseli} received the B.Sc. and M.Sc. degrees in
electrical engineering from the Amirkabir University of Technology,
Tehran, Iran, in 2013 and 2015, resp. He also received the
Ph.D. degree in Engineering Sciences (Mechanical Engineering) from UC
San Diego, USA, in 2022, where he is currently a postdoctoral
researcher. His research interests include system identification,
nonlinear systems, network systems, data-driven modeling and control,
and distributed and parallel computing. Dr. Haseli is the recipient of
the Bronze Medal of the 2014 Iran National Mathematics Competition and
the Best Student Paper Award of the 2021 American Control Conference.

\smallskip
\begin{wrapfigure}{l}{0.25\linewidth}
\centering
\includegraphics[width=\linewidth]{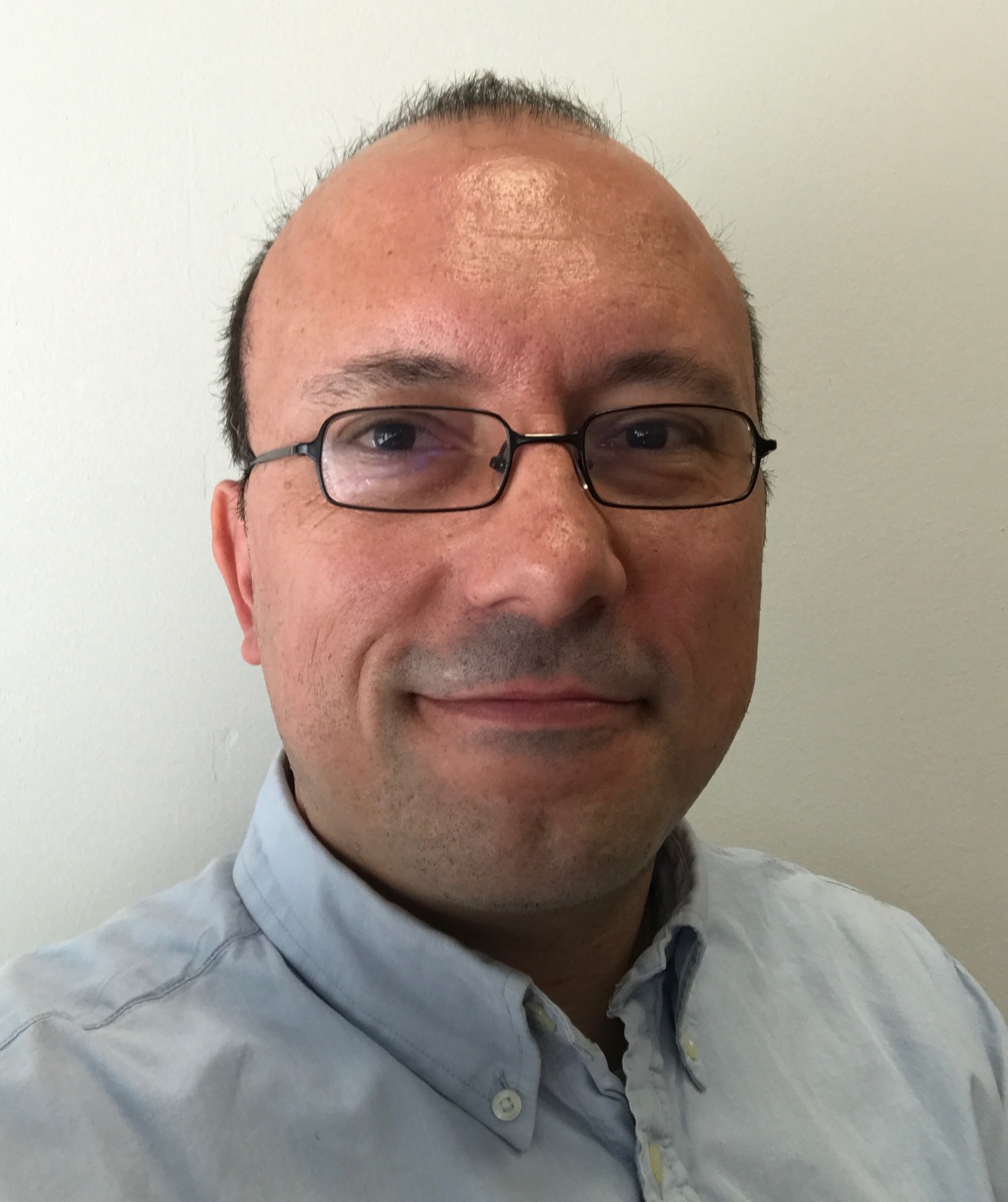}
\end{wrapfigure}
\textbf{Jorge Cort\'es} received the Licenciatura degree in
mathematics from Universidad de Zaragoza, Spain, in 1997, and the
Ph.D. degree in engineering mathematics from Universidad Carlos III de
Madrid, Spain, in 2001. He held postdoctoral positions with the
University of Twente, The Netherlands, and the University of Illinois
at Urbana-Champaign, USA.
He is a Professor in the Department of Mechanical and Aerospace
Engineering, UC San Diego, USA. He is a Fellow of IEEE, SIAM, and
IFAC. His research interests include distributed control and
optimization, network science, nonsmooth analysis, reasoning and
decision making under uncertainty, network neuroscience, and
multi-agent coordination in robotic, power, and transportation
networks.

\end{document}